\documentclass{IEEEojcsys}
\UseRawInputEncoding
\usepackage[colorlinks,urlcolor=blue,linkcolor=blue,citecolor=blue]{hyperref}

\usepackage{amsmath, amssymb, amsfonts, amsthm}
\usepackage{mathtools}
\usepackage{bbold}
\usepackage{textcomp}

\usepackage{graphicx}
\usepackage{epsfig}          
\usepackage{grffile}
\usepackage{svg}
\usepackage{psfrag}
\usepackage{float}
\usepackage{stfloats}
\usepackage{subcaption}
\usepackage{caption}
\usepackage{tikz}
\usepackage{pgfplots}
\usepackage{circuitikz}
\pgfplotsset{compat=newest}
\usepgfplotslibrary{patchplots}

\usetikzlibrary{
    shapes,
    shapes.geometric,
    arrows,
    positioning,
    fit,
    plotmarks,
    calc
}

\usepackage{tabularx}
\usepackage{booktabs}
\usepackage{multirow}
\usepackage{array}
\newcolumntype{Y}{>{\centering\arraybackslash}X}

\usepackage{algorithm}
\usepackage{algpseudocode}
\usepackage{listings}

\usepackage{cite}
\usepackage{enumerate}
\usepackage{multicol}
\usepackage{color}

\usepackage{hyperref}

\newtheorem{problem}{Problem}
\newtheorem{definition}{Definition}
\newtheorem{lemma}{Lemma}
\newtheorem{remark}{Remark}
\usepackage{xcolor} 

\definecolor{myblue}{RGB}{0, 114, 178}
\definecolor{myred}{RGB}{214, 39, 40}
\definecolor{myorange}{RGB}{240, 140, 0}
\definecolor{mygreen}{RGB}{44, 160, 44}   
\definecolor{mypurple}{RGB}{148, 103, 189} 
\definecolor{myblack}{RGB}{0, 0, 0} 

\begin{document}

\title{Mamba Sequence Modeling meets \\ Model Predictive Control} 


\author{M. Cevaal\affilmark{1}} 

\author{T.O. de Jong\affilmark{2} (Student Member, IEEE)}

\author{M. Lazar\affilmark{2} (Senior Member, IEEE)}

\affil{Tilburg University, The Netherlands} 
\affil{Eindhoven University of Technology, The Netherlands}

\corresp{CORRESPONDING AUTHOR: M. Cevaal (e-mail: \href{mailto:m.cevaal@student.tue.nl}{m.cevaal@tilburguniversity.edu})}
\authornote{}

\markboth{Mamba Sequence Modeling meets Model Predictive Control}{M.Cevaal {\itshape ET AL}.}

\begin{abstract}
In this paper, we consider the design of Model Predictive Control (MPC) algorithms based on Mamba neural networks. Mamba is a neural network architecture capable of sub-quadratic computational scaling in sequence length with state-of-the-art modeling capabilities.
We provide a consistent and complete mathematical description of the Mamba neural network.
Then, adjustments and optimizations are made to construct a decoder-only Mamba multi-step predictor for MPC and an input-output formulation is given for sequence-to-sequence modeling of dynamical systems. The performance of Mamba-MPC is evaluated on several numerical examples and compared to a Long-Short-Term-Memory-Based MPC (LSTM-MPC) and a Transformer based MPC (Transformer-MPC) for computational scaling comparison. First, a Single-Input-Single-Output (SISO) Van der Pol oscillator is considered, where stability, reference tracking, and noise robustness are evaluated. Then, a MIMO Four Tank setup is introduced where Multiple-Input-Multiple-Output (MIMO) reference tracking is evaluated. Lastly, Mamba-MPC is implemented on a physical Quanser Aero2 setup for closed-loop reference tracking. The results demonstrate that Mamba-MPC is able to stabilize and track a reference for SISO and MIMO systems, both in simulation and on a physical setup. Moreover, Mamba-MPC consistently outperforms LSTM-MPC in predictive control and is significantly computationally faster.
\end{abstract}


\maketitle

\section{INTRODUCTION}
Model Predictive Control (MPC) is a highly successful advanced control methodology which can be merged with \emph{Neural Networks} (NN) as predictors and thus leverage their universal approximation properties \cite{Chen1995} for closed-loop control. Particularly when operating on complex nonlinear systems, this is an interesting approach because accurate first-principles models are often difficult to obtain.
 Nevertheless, a crucial limiting factor for the successful application of NN within the MPC framework is that their intricate nonlinear architecture inherently introduces computational complexity, which restricts computational tractability and poses challenges for online implementation.
Early efforts to model nonlinear dynamics for MPC used feedforward NNs \cite{Lazar2002}, but these models struggle to capture long-range temporal dependencies and suffer from cascading errors in multi-step predictions. Recurrent architectures such as \emph{Recurrent Neural Networks} (RNN) for MPC \cite{Pan2009} address this by propagating temporal information via latent states and have been shown to explicitly exhibit state-space dynamics in \cite{Hoekstra2023}. Nevertheless, modeling long-range dependencies with RNN remains a challenge. Consequently, gating augmentations were introduced in \emph{Long Short-Term Memory} (LSTM) \cite{Hochreiter1997, Mohajerin2018} and \emph{Gated Recurrent Units} (GRU) \cite{Cho2014, Bonassi2021} to enable long-range prediction. Despite these advances, accurate modeling over longer prediction horizons remains challenging due to the persistent accumulation of errors inherent to one-step ahead predictors.
More recently, Transformer networks \cite{Vaswani2017} have emerged for modeling global temporal relations which is particularly advantageous in the modeling of extremely long prediction horizons. A multi-step Transformer predictor for MPC was constructed in \cite{Park2023}, and compared to LSTM. In comparison to LSTM, the multi-step Transformer showed superior performance and speed, although the computational benefit largely stemmed from the construction of a  multi-step predictor. A challenge with Transformers is that the main mechanism of Transformers, attention, leads to quadratic computational scaling when increasing the length of the prediction horizon. This limits the practical implementation of Transformers in MPC.
A novel NN architecture, Mamba \cite{Gu2023}, has recently emerged as a compelling alternative, built on structured state-space models \cite{Gu2021}.
Mamba incorporates a context-aware selection mechanism that propagates relevant information from past sequence steps in the hidden states, effectively incorporating the gating mechanism found in GRUs. Gating mechanisms in NNs have been shown to enable the network to learn time warping functions \cite{Tallec2018}, which aids efficient storing of information. Crucially, the selection mechanism enables Mamba to scale linearly in computational cost when increasing the sequence length. Furthermore, Mamba has demonstrated superior performance compared to Transformers and DeepONet \cite{Lu2020} in the modeling of dynamical systems \cite{Hu2024} when used as a neural operator for learning \emph{Sequence-to-Sequence} (Seq2Seq) mappings. Consequently, a Mamba-based prediction model is an excellent candidate for predictive control, which, to the best of our knowledge, has not been explored. 

It is worth mentioning that an approach that combines elements of the Mamba architecture with a linear-in-control Koopman operator approach for MPC was recently presented in \cite{Li-Zhaoyang}. Therein Mamba is used as a linear state-space matrix generator in combination with a Koopman operator to form a lifted time-varying linear state-space system model. 
To the best of our knowledge, an implementation of Mamba for MPC where the full Mamba architecture is employed as a multi-step predictor has not yet been explored.

Therefore, in this paper, we investigate the application of Mamba in predictive control of nonlinear dynamical systems by constructing a data-driven Mamba-MPC algorithm. The first challenge is to construct a Mamba-based multi-step predictor while maintaining equal sequence length between input and output. 
The second challenge is to construct a multi-step formulation of Seq2Seq to enable input/output mapping for a set of initial conditions. We thus aim to formulate a Mamba-based predictor that allows for MIMO multi-step predictions. 
The first contribution of this paper is \emph{(i)} a complete and consistent mathematical description of the underlying operations employed within the Mamba NN architecture. Note that, to the best of our knowledge, a complete rigorous mathematical description of the Mamba NN is not available in the literature, which currently hampers understanding of the underlying functional operations and reproducibility. Furthermore, we include a Lipschitz boundedness analysis of the complete Mamba neural network architecture. The second contribution is \emph{(ii)} to provide an adaptation and optimization of Mamba for Mamba-MPC as well as a Seq2Seq input-output embedding. The third contribution is \emph{(iii)} implementation of Mamba-MPC in real-time control of a challenging physical Quanser Aero2 system. To enhance reproducibility, we provide for all learning and MPC algorithms the implementation code in Python using PyTorch and CasADI. 

The remainder of the paper is structured as follows. Section \ref{sec:Notations} introduces notations and definitions. Section \ref{sec:Prelim} reviews MPC preliminaries and define the problem statement. In Section \ref{sec:Mamba} a mathematical description of Mamba is provided. In Section \ref{sec:Seq2Seq}, the Mamba-MPC algorithm is derived. In Section \ref{sec:NumExample}, we test Mamba-MPC for 2 numerical examples, the Van der Pol oscillator and a MIMO Four Tank system. In Section \ref{sec:AERO} we provide an implementation on a physical setup as a proof of concept for practical application of Mamba-MPC. Section \ref{sec:Conclusion} presents conclusions and future work.

\section{Notations and Definitions}\label{sec:Notations}
Basic notation and definitions are given next.
\paragraph{Vectors}
For a finite collection of vectors $\xi_i \in \mathbb{R}^{n_i}$ ($i=1, \dots, q$), let $\operatorname{col}(\xi_1, \dots, \xi_q) := [\xi_1^\top \ \dots \ \xi_q^\top]^\top$.

\paragraph{Matrices}
For any two matrices $A \in \mathbb{R}^{n \times m}$ $B \in \mathbb{R}^{n \times m}$, we denote the Hadamard product by $A \odot B$, defined as 
\begin{equation*}
    A \odot B = \begin{bmatrix}
        a_{1,1}b_{1,1}  & \cdots & a_{1,m}b_{1,m} \\
        \vdots & & \vdots \\
        a_{n,1}b_{n,1} & \cdots & a_{n,m}b_{n,m}          
    \end{bmatrix}
\end{equation*}
which can be extended to high-order tensors.
Given scalars $a_1, a_2, \dots, a_n\in \mathbb{R}$ we define the diagonal matrix $A = \text{diag}\left( a_1, a_2, \dots, a_n\right)\in \mathbb{R}^{n\times n}$ where $A_{i,j} = a_i$ when $i=j$ and $A_{i,j} = 0$ when $i \neq j$ such that $A$ has $a_1, a_2, \dots, a_n$ on its diagonal and zeros elsewhere.
Given a vector $w = \left[w_{1}, \dots, w_{l}\right]^{\top} \in \mathbb{R}^{l}$ we define a Toeplitz matrix constructed from vector $w$ as
\begin{equation}
\text{Toeplitz}(w) = 
    \begin{bmatrix}
        w_{1}  & \dots& w_{l} & 0 & \dots & 0 \\
        0 &    w_{1}  & \dots& w_{l} & \dots & 0 \\
        \vdots & \vdots & \ddots & \ddots & \ddots &  \vdots \\
        0 & \dots & 0 &  w_{1}  & \dots& w_{l}
    \end{bmatrix}.
\nonumber
\end{equation}
\paragraph{Tensors}
In this paper tensors up to rank $3$ are mostly considered. Let calligraphic variables, e.g.,  $\mathcal{T} \in \mathbb{R}^{d_1 \times d_{2} \times d_{3}}$ denote a tensor. When operating on a tensor, let $\mathcal{T}_{:,i, :}\in \mathbb{R}^{d_1 \times  d_{3}}$, denote the $i^{th}$ slice of the second dimension, e.g., for a tensor $\mathcal{B}\in \mathbb{R}^{2\times 2\times 2}$, e.g.,
\begin{equation}
    \mathcal{B} = \left(\left\{ \begin{bmatrix}
        b_1 & b_2 \\
        b_3 & b_4
    \end{bmatrix} ,
    \begin{bmatrix}
        b_5 & b_6 \\
        b_7 & b_8 
    \end{bmatrix} \right\}\right),\nonumber
\end{equation}
we denote the $1$-st slice of dimension $d_2$ by:
\begin{equation}
  \mathcal{B}_{:,1,:}
    = \left\{\begin{bmatrix}
        b_1 & b_2 
    \end{bmatrix},
        \begin{bmatrix}
        b_5 & b_6 
    \end{bmatrix}\right\} \nonumber.
\end{equation}
Furthermore, we denote the flattening of a rank 3 tensor by $\text{mat}:\mathbb{R}^{b\times n \times m}\rightarrow \mathbb{R}^{bn \times m}$, which amounts to combining the first and second dimension. For example,
\begin{equation}
    \text{mat}\left(\left\{ \begin{bmatrix}
        a_1 & a_2 \\
        a_3 & a_4
    \end{bmatrix} ,
    \begin{bmatrix}
        a_5 & a_6 \\
        a_7 & a_8 
    \end{bmatrix} \right\}\right) 
    = \begin{bmatrix}
                a_1 & a_2 \\
        a_3 & a_4 \\
                a_5 & a_6 \\
        a_7 & a_8 
    \end{bmatrix}\nonumber.
\end{equation}
For a tensor  composed of matrices $ \mathcal{Q}= \left(Q_1, Q_2, \dots, Q_M\right)$, the construction of a block diagonal matrix from tensor $\mathcal{Q}$ is denoted by
\begin{equation}
\text{blk}\left(\mathcal{Q}\right) = \begin{bmatrix}
    Q_1 & 0 & \dots & 0 \\
    0 & Q_2 & \dots & 0 \\
    \vdots  & \vdots & \ddots& \vdots \\
    0& \dots & 0 &  Q_M
\end{bmatrix}.
\nonumber
\end{equation}
Given a matrix $F\in \mathbb{R}^{m\times n}$ where $F = \left[f_1,f_2,\dots, f_n\right]$ and $f_i$ is the $i$-th column of $F$, the vectorization operation on $F$ is defined as
\begin{equation}
    \text{vec}\left(F\right) = \left[f_1^{\top}, f_2^{\top}, \dots, f_n^{\top}\right]^{\top}.
    \nonumber
\end{equation}
Moreover, for the  process where a matrix is created from a vector $g = \left[g_1, \dots, g_l \right] \in \mathbb{R}^{mn}$ is defined as
\begin{equation}
    H =  Q_{V\to M}\left(g,  m, n\right),
    \nonumber
\end{equation}
where $(m, n)$ denote the dimension of $H$, i.e., $H\in \mathbb{R}^{m\times n}$. 

\paragraph{Broadcasting Hadamard Product}
In order to perform the Hadamard product between tensors of different dimensions we define the broadcasting operation and the broadcasting Hadamard product. Consider two rank $3$ tensors $\mathcal{D} \in \mathbb{R}^{L_1 \times L_2 \times L_3}$ and $\mathcal{E} \in \mathbb{R}^{J_1 \times J_2  \times J_3}$. Then the broadcasting operation can be used to perform the Hadamard product between tensors of different dimensions, $\mathcal{F}=\mathcal{D} \boxdot \mathcal{E}$,  as follows, 
\begin{align}
&\mathcal{F}= \mathcal{D} \boxdot\mathcal{E} := \text{bc}\left(\mathcal{D},\text{size}\left(\mathcal{E}\right)\right)  \odot \text{bc}\left(\mathcal{E},\text{size}\left(\mathcal{D}\right)\right), \nonumber \\
& \mathcal{F} \in \mathbb{R}^{\max \left(L_1, J_1\right) \times \max \left(L_2, J_2\right) \times \max \left(L_3, J_3\right)} .\nonumber 
\end{align}   
Above $\text{size}: \mathbb{R}^{n\times m \times p}\rightarrow \mathbb{N}\times  \mathbb{N}\times \mathbb{N}$ returns the shape of the input tensor as an $n$-tuple such that $\text{size}\left(\mathcal{J}\right) = \left(J_{1},J_{2}, J_{3}\right)$. Furthermore, the broadcasting function $\text{bc}:\mathbb{R}^{L_1\times L_2 \times L_3} \times (\mathbb{N} \times \mathbb{N}  \times \mathbb{N} ) \rightarrow \mathbb{R}^{K_1\times K_2 \times K_3}$, where $K_1\geq L_1$, $K_2\geq L_2$ and $K_3\geq L_3$, is a function that has a tensor and $n$-tuple as arguments such that for tensor $\mathcal{D}$ and $\text{size}\left(\mathcal{E}\right)$ for all $n$, if $L_{n}=1$, replicates all elements of $\mathcal{D}$ along the $n$-th mode $J_{n}$ times such that $\mathcal{D}' = \text{bc}\left(\mathcal{D}, \text{size}\left(\mathcal{E}\right)\right)$.
To illustrate the application of the broadcasting Hadamard product, consider  the following example for a matrix and a column vector:
\begin{equation}
    \begin{bmatrix}
    f_1 & f_2 \\
    f_3  & f_4
\end{bmatrix} \boxdot \begin{bmatrix}
    g_1 \\
    g_2
\end{bmatrix} = 
    \begin{bmatrix}
    f_1 & f_2 \\
    f_3  & f_4
\end{bmatrix} \odot \begin{bmatrix}
    g_1 & g_1 \\
    g_2 & g_2
\end{bmatrix}\nonumber.
\end{equation}
Similarly, for a matrix and row vector:
\begin{equation}
    \begin{bmatrix}
    f_1 & f_2 \\
    f_3  & f_4
\end{bmatrix} \boxdot \begin{bmatrix}
    g_1 & g_2 \\
\end{bmatrix} = 
    \begin{bmatrix}
    f_1 & f_2 \\
    f_3  & f_4
\end{bmatrix} \odot \begin{bmatrix}
    g_1 & g_2 \\
    g_1 & g_2
\end{bmatrix} \nonumber.
\end{equation}
For a tensor operation consider the following example,
\begin{align}
\left\{ \begin{bmatrix}
        a_1  \\
        a_2
    \end{bmatrix} ,
    \begin{bmatrix}
        a_3 \\
        a_4 
    \end{bmatrix} \right\}
        \boxdot
    \left\{ \begin{bmatrix}
        a_5 & a_5 \\
        a_6 & a_6
    \end{bmatrix} ,
    \begin{bmatrix}
        a_7 & a_7\\
        a_8 & a_8
    \end{bmatrix}\right\} = \nonumber \\
\left\{ \begin{bmatrix}
        a_1 & a_1  \\
        a_2 & a_2
    \end{bmatrix} ,
    \begin{bmatrix}
        a_3 & a_3 \\
        a_4 & a_4
    \end{bmatrix} \right\}
        \odot
\left\{ \begin{bmatrix}
        a_5 & a_5 \\
        a_6 & a_6
    \end{bmatrix} ,
    \begin{bmatrix}
        a_7 & a_7\\
        a_8 & a_8
    \end{bmatrix}\right\} \nonumber.
\end{align}
For a formal definition consider the following work \cite{matsui2026}. To perform the broadcasted Hadamard product, a necessary condition must be satisfied for the broadcasted Hadamard product to be well-defined.  The broadcasting condition for $\mathcal{D}$ and $\mathcal{E}$ holds if for any $n\in \left\{1,2,3 \right\}$ one of the following conditions is satisfied: $L_{n} = J_{n}$, $L_{n}=1$ or $J_{n}=1$.

\paragraph{Einstein Summation}
To perform Einstein summation, consider two tensors $\mathcal{H} \in \mathbb{R}^{P \times T \times R \times E}$ and $K \in \mathbb{R}^{P \times T \times E}$. The output $L \in \mathbb{R}^{P \times T \times R}$ is obtained by contracting the common dimension $E$ as follows:
\begin{equation}
L_{p,t,r} = \sum_{e=1}^{E} \mathcal{H}_{p,t,r,e} K_{p,t,e},\nonumber
\end{equation}
where the indices $p, t, r$ correspond to the free dimensions and $e$ is the summation index.

\section{MPC Preliminaries and Problem Statement}\label{sec:Prelim}
Consider a continuous time-invariant nonlinear MIMO state-space system, i.e., 
\begin{subequations}\label{ContinuousNLMIMO}
 \begin{align}
    \dot{x}\left(t\right) &= f_c\left(x\left(t\right), u\left(t\right)\right) \\
    y\left(t\right) &= h_c\left(x\left(t\right), u\left(t\right) \right)
\end{align}
\end{subequations}
with inputs $u\in \mathbb{U }\subset \mathbb{R}^{n_u}$, states $x \in \mathbb{X}\subset \mathbb{R}^{n_x}$ and outputs $y \in \mathbb{Y}\subset \mathbb{R}^{n_y}$, time $t \in \mathbb{R}_{\geq 0}$. Above $f_c: \mathbb{U}\times \mathbb{X}\to \mathbb{Y}$ and $h_c: \mathbb{U}\times \mathbb{X}\to \mathbb{Y}$ are unknown suitable functions that allow for existence of unique solutions. 

To perform sequential computation on \eqref{ContinuousNLMIMO}, the model needs to be discretized, for example, by means of zero-order hold (ZOH). 
In general, a  discrete-time nonlinear MIMO system is described by, i.e.,
\begin{subequations}    \label{eq:NLMIMO}
 \begin{align}
    x\left(k+1\right) &= f\left(x\left(k\right), u\left(k\right)\right), \\
    y\left(k\right) &= h\left(x\left(k\right), u\left(k\right) \right),
\end{align}   
\end{subequations}
where for a discrete time instant $k \in \mathbb{N}$, the continuous time $t=kT_s$, where $T_s \in \mathbb{R}$ represents the sampling time/period.
In discrete-time predictive control, sequences of future inputs and predicted outputs over a prediction horizon $N \in \mathbb{N}_{\geq 1}$ are generally of interest.
\begin{figure}[b]
\begin{tikzpicture}[
Block1/.style={draw=blue, fill = blue!5, very thick,  minimum width=2cm, minimum height=1.5cm},
Block2/.style={draw=red, fill = red!5, very thick,  minimum width=2cm, minimum height=1.5cm},
Block3/.style={draw=red, fill = red!5, very thick,  minimum width=1.2cm, minimum height=1.2cm},
Block4/.style={draw=black, fill = black!5, very thick,  minimum width=1.5cm, minimum height=1.5cm},
node distance=10mm  
]
\node[Block2] (block3) {\parbox{3cm}{\centering Model Predictive Control \\ (Problem~\ref{eq:PrototypeMPC})}};
\node[Block2] (block4) [below=of block3] {\parbox{3cm}{\centering \begin{align*}
    \dot x &= f_c\left(x,u\right) \\
    y &= h_c\left(x\right)
\end{align*}}};
 \node[Block3] (block1_small) [right=of $(block3)!0.5!(block4)$, xshift=15mm, inner sep=0pt, outer sep=0pt] {ZOH};
\node[Block3, inner sep=0pt, outer sep=0pt, minimum width=1.2cm, minimum height=1.2cm] (block2_small) [left=of $(block3)!0.5!(block4)$, xshift=-15mm] { 
   \begin{circuitikz}
     \draw (0.1,0.2) to[nos, o-o] (1.0,0.2);  
     \node[above=0.0cm] at (0.5, 0.1) {$T_s$};  
   \end{circuitikz}
}; 
\draw[->, very thick] (block4.west) -| (block2_small.south) node[midway, right, xshift=2.5mm, yshift=3mm] {$y(t)$};
 
\draw[->, very thick] (block2_small.north) |- (block3.west) node[midway, right, xshift=0mm, yshift=3mm] {$y(kT_s)$};
 \draw[->, very thick] (block3.east) -| (block1_small.north) node[midway, right, xshift=-10.0mm, yshift = 3mm] {$u(kT_s)$};
 \draw[->, very thick] (block1_small.south) |- (block4.east) node[midway, left, xshift= -1.0mm, yshift = 3.0mm] {$u(t)$};
 \draw[->, very thick] 
  ([yshift=25pt]block3.north) -- (block3.north) 
  node[midway, right] 
  {$\mathbb{F}\left(u_{\left[0,N-1\right]}\left(k\right),  x_{\mathrm{0}}\left(k\right)\right)$}; 
\end{tikzpicture}
    \caption{Illustration of the Model Predictive Control closed-loop control architecture.}
    \label{fig:MPCOverview}
\end{figure}
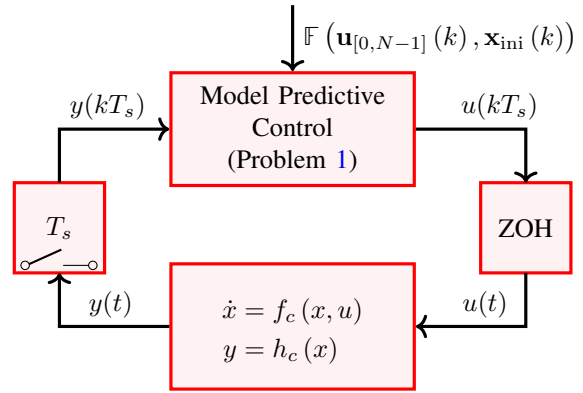

Consider thus a sequence of inputs $u_{\left[0,N-1\right]}\left(k\right)$ and outputs $y_{\left[1,N\right]}\left(k\right)$ described by:
\begin{subequations}
    \begin{align}
        u_{\left[0,N-1\right]}\left(k\right) &:= \operatorname{col}(\left[u\left(0|k\right), \dots,u\left(N-1|k\right) \right]), \\
y_{\left[1,N\right]}(k) &:=\operatorname{col}(\left[y\left(1|k\right), \dots,y\left(N|k\right) \right]),
    \end{align}
\end{subequations}
where the input sequence is mapped to a sequence of outputs for any initial condition by a suitable map $\mathbb{F}: \mathbb{R}^{n_u N}\times \mathbb{R}^{n_x}\to \mathbb{R}^{n_yN}$.
Assuming that a suitable map $\mathbb{F}$ can accurately describe the Seq2Seq multi-step representation of the MIMO nonlinear system in \eqref{eq:NLMIMO}, then a prototype predictive control problem can be formulated as  follows: 
\begin{problem}{Prototype MPC problem} \label{eq:PrototypeMPC}
\begin{equation}
\label{prob:NOMPC}
\begin{split}
&\min _{u_{\left[0,N-1\right]}\left(k\right), y_{\left[1,N\right]}\left(k\right)} l_N\left(y\left(N|k\right)\right)+\sum_{i=0}^{N-1} l\left(y\left(i|k\right), u\left(i|k\right)\right) \\
&\text{ subject to: }\\
& y_{\left[1,N\right]}(k)= \mathbb{F}\left(u_{\left[0,N-1\right]}\left(k\right),  x_{0}\left(k\right)\right),\\
&u_{\left[0,N-1\right]}\left(k\right)\in\mathbb{U}\times\ldots\times\mathbb{U},\\
&y_{\left[1,N\right]}\left(k\right)\in\mathbb{Y}\times\ldots\times\mathbb{Y},\\
&y(0|k)=y(k).
\end{split}
\end{equation}    
\end{problem}
In \eqref{prob:NOMPC}, $x_{0}\left(k\right)$ is the initial condition that can be the initial state $x\left(k\right)$, if the measured states are available. $\mathbb{F}$ is a suitable map such that the sequence of inputs $u_{\left[0,N-1\right]}\left(k\right)$ and the initial condition $x_{0}(k)$ are mapped to the output sequence $y_{\left[1,N\right]}\left(k\right)$. $\mathbb{F}$ can be explicitly constructed if the functions $f,h$ are known, by repeated appropriate function compositions. An illustration of using Problem~\ref{eq:PrototypeMPC} to define a closed-loop control architecture is provided in Figure~\ref{fig:MPCOverview}.

For the learning of the map $\mathbb{F}$, several approaches can be considered. A common data-driven approach is to utilize NN architectures for the multi-step predictor. A feedforward NN was applied as multi-step predictor in MPC \cite{Lazar2002} and achieved promising results. 
In \cite{Park2023} a multi-step predictor is constructed based on a Transformer network and is successfully applied in the MPC framework but at significant computational cost.
Another approach is operator networks such as DeepONet \cite{Lu2020}, which were recently adapted for implementation in MPC in \cite{deJong2025}. 
In \cite{Hu2024} a comparative study is done involving DeepONet, Transformer based approaches, and Mamba for the modeling of dynamical systems. Mamba consistently achieves superior performance with respect to DeepONet and Transformer based approaches. Furthermore, Mamba was inherently designed for linear computational scaling with respect to the prediction horizon. Consequently, Mamba is of significant interest as an application to learn nonlinear mappings for multi-step prediction in the MPC framework. 
In particular, Mamba-based decoder-only architectures are of interest, such that the total parameter count and therefore computational complexity in Problem~\ref{eq:PrototypeMPC} remains minimal.
Therefore, the problem statement can be summarized as follows. 

\textit{Problem statement:}\quad Given the discrete-time nonlinear MIMO system in \eqref{eq:NLMIMO}, the goal is to 
approximate the multi-step Seq2Seq map $\mathbb{F}$ from future inputs $u_{\left[0,N-1\right]}(k)$ to future outputs $y_{\left[1,N\right]}\left(k\right)$ for a set of initial conditions $x_{\mathrm{0}}\left(k\right)$  by means of Mamba NNs. The input of the network should have a sequence length equal to the sequence length of the output such that decoder-only structures can be utilized for computational tractability. The model should predict a multi-step sequence in a single forward pass, such that efficient predictive control is enabled.   

\begin{remark}
    If state data is unavailable, a a stacked vector of past inputs and outputs can be used that encode the state information. In the rest of the paper we assume that state data is available.
\end{remark}

\begin{figure}
    \centering
    \includegraphics[width=1\linewidth]{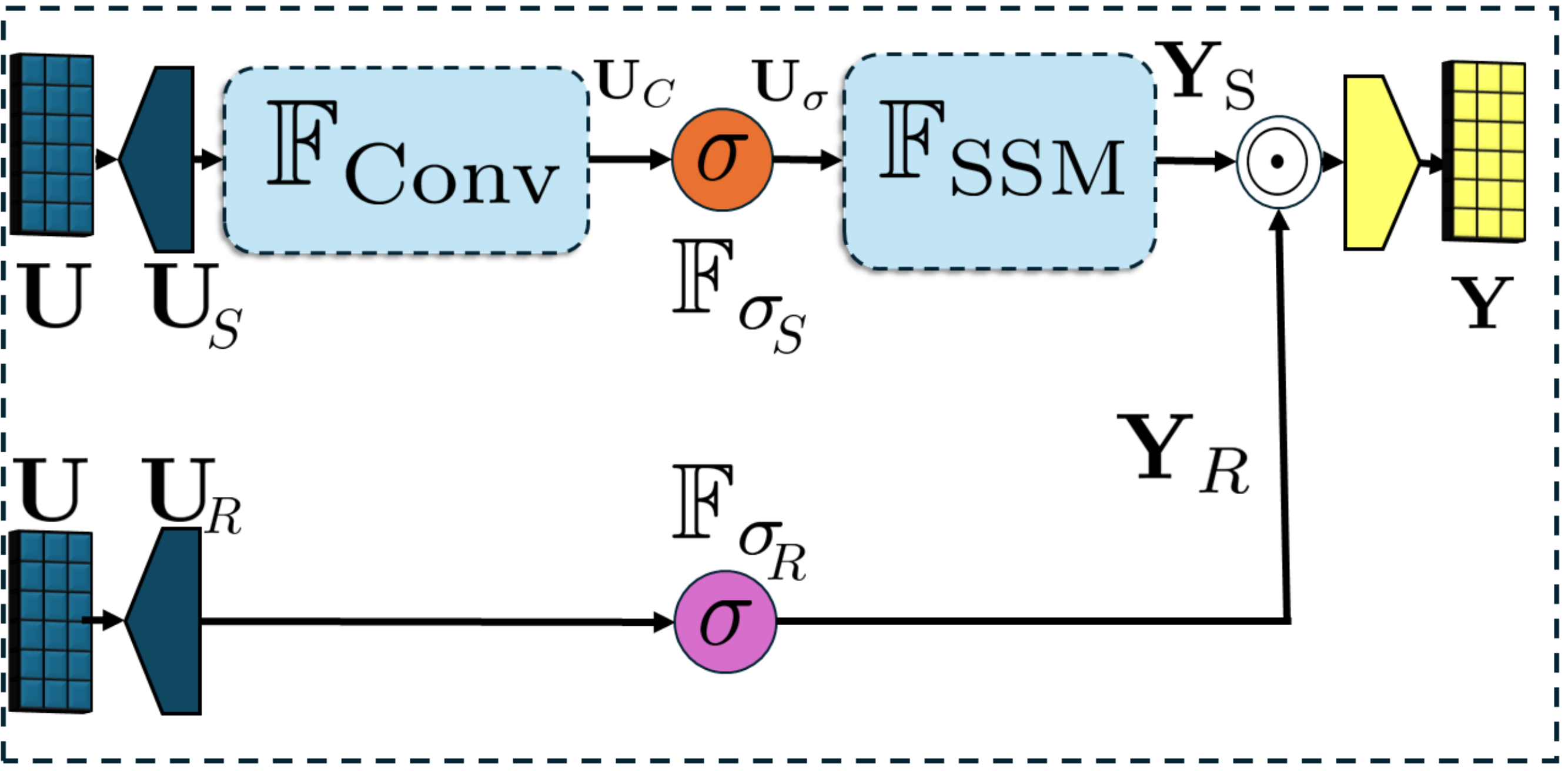}
    \caption{A schematic of the Mamba block architecture. The input $U$ is processed through two parallel paths: a main path featuring a convolution layer, a hidden layer of activation functions, a SSM block that implements a time-varying linear dynamical system in the feature space, and a residual path with a hidden layer of activation functions. The final output is the element-wise product of the outputs of these two paths.}
    \label{fig:MambaFunctions}
\end{figure}

\section{Mamba: Mathematical Description}\label{sec:Mamba}
Let the map from input to output be described by $Y = \mathbb{F}_{\mathrm{Mamba}}\left(U\right)$ where  $\mathbb{F}_{\mathrm{Mamba}}\left(U\right)$ is defined according to:
\begin{equation}
        \mathbb{F}_{\mathrm{Mamba}}\left(U\right):=  \left(\mathbb{F}_{SSM} \circ \sigma\left(\mathbb{F}_{\mathrm{Conv}}\left(U W_S^{\top} \right)\right) \odot Y_R\right)W_Y^{\top} \label{eq:Fmamba} 
\end{equation}
where
\begin{subequations}
    \begin{align}
        Y_R &= \sigma\left(UW_R^{\top}\right), \label{eq:Residual} \\
        \sigma(U_C) &= U_{C}\odot\frac{1}{1+\exp\left(-U_{C}\right)}. \label{eq:sigma}
    \end{align}
\end{subequations}
In $F_{\mathrm{Mamba}}$, the functions $\sigma: \mathbb{R} \to  \mathbb{R}$ in \eqref{eq:Fmamba} and \eqref{eq:Residual} refer to an element-wise \emph{Sigmoid Linear Unit} (SiLU) activation function\cite{Shazeer2020}. Consequently, the exponential function in \eqref{eq:sigma} is an element-wise function.
Furthermore, mappings $W_S^{\top}, W_R^{\top}: \mathbb{R}^{L \times D} \to  \mathbb{R}^{L \times ED}$ are lifting operators that expand the feature space of the input, while $W_Y^{\top}: \mathbb{R}^{L \times ED} \to  \mathbb{R}^{L \times D}$ are projection operators that reduce the feature space.
The functions $\mathbb{F}_{SSM}$ and $ \mathbb{F}_{\mathrm{Conv}}$ require a more detailed description and will be further highlighted in the following sections. 
To illustrate the general structure of the Mamba function, a graphical description can be seen in Figure \ref{fig:MambaFunctions}.
Within the architecture of Mamba, an input is typically provided as a tensor denoted by $\mathcal{U} \in \mathbb{R}^{B_{a} \times L \times D}$ where the parameter definitions are given in Table \ref{tab:Definitions}. For the following section, the dimension $B_{a}$ will be considered of size one, such that $U = \text{mat}(\mathcal{U})\in \mathbb{R}^{L \times D}$, for simplicity.

\begin{table}[t!]
\centering
\caption{Model Parameter Definitions}
\begin{tabular}{|l|l|}
\hline
\textbf{Parameter} & \textbf{Definition} \\
\hline
$B_{a}$ & Batch Size \\
$L$ & Sequence Length \\
$D$ & Model Dimension \\
$E$ & Expansion  Factor \\
$N$ & Prediction Horizon \\
$S$ & Hidden State Dimension \\
$K$ & Convolutional Kernel Size \\
$d_{\tau}$ & Rank of $\Delta$ Approximation\\
\hline
\end{tabular}
\label{tab:Definitions}
\end{table}

\subsection{1D Convolution}
\begin{figure}
    \centering
    \includegraphics[width=1\linewidth]{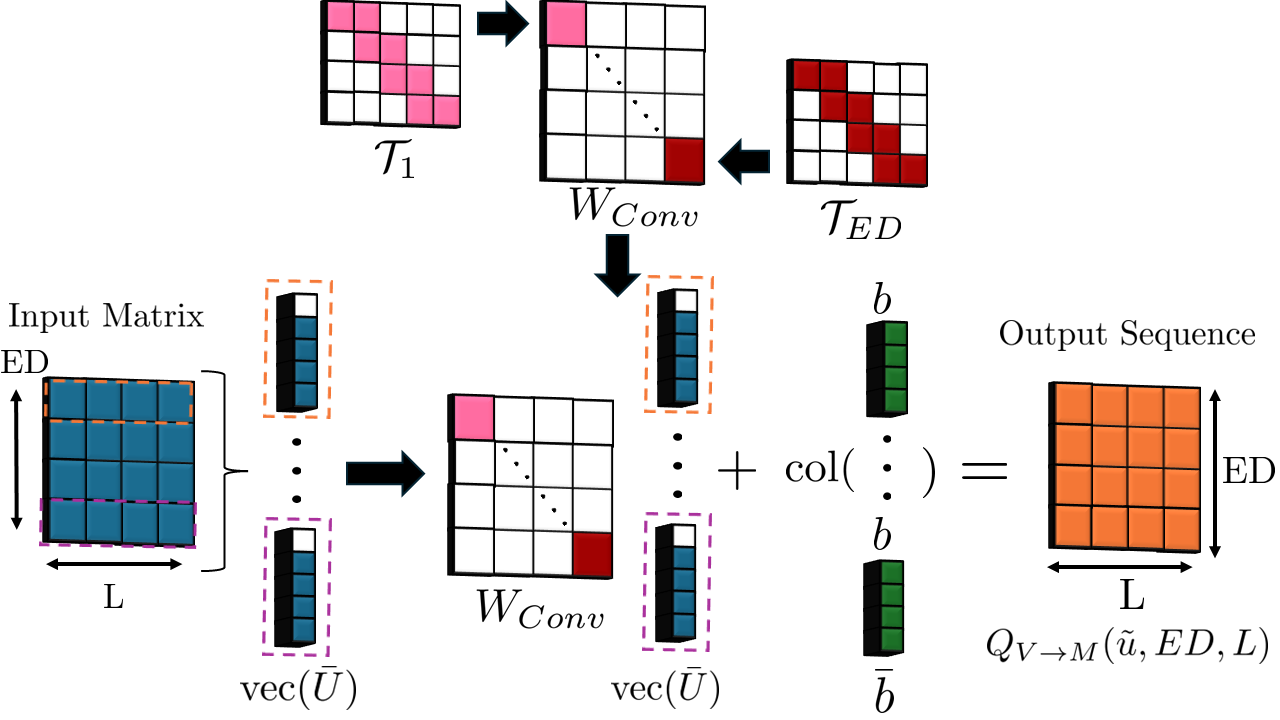}
    \caption{An illustration of the 1D depthwise causal convolution operation implemented as a matrix multiplication. The convolution kernel for each dimension $d$ is transformed into a Toeplitz matrix $T_d$, which are then combined into a block-diagonal matrix $T$. This allows the entire depthwise convolution to be computed efficiently in parallel.}
    \label{fig:Conv1D}
\end{figure}

Next, a causal depthwise 1 dimensional convolutional layer (1D Conv) is applied along the sequence length dimension $L$ of the expanded input. The addition of 1D Conv in Mamba is such that linear combinations of features of prior sequence steps can be incorporated into the current sequence step. Hence, this aids in the modeling linear temporal dependencies.
An illustration of the operation is given in Figure~\ref{fig:Conv1D}. 

Let the 1D Conv function $ \mathbb{F}_{\mathrm{Conv}}$ be described by:
\begin{equation}
    \mathbb{F}_{\mathrm{Conv}}\left(U_{S}\right) := Q_{V\to M}\left(\tilde{u}, ED, L\right)^\top,
    \label{eq:Conv1D}
\end{equation}
where
\begin{subequations}
\begin{align}
    \tilde{u}=& W_{Conv}\text{vec}\left(\tilde{U}\right)+\tilde{b}, \\
        W_{Conv} =& \text{blk}\left(\mathcal{T} \right), \\
    \mathcal{T} =& \text{Toeplitz}\left( \kappa_{d,:}\right), \,  d = 1,\dots, ED, \\
\tilde{U} =& \left[\bar{0}_{l}, U^{\top}_{S}, \bar{0}_r\right], \\
\tilde{b} =& \text{col}\left(b, \dots, b\right) \in \mathbb{R}^{ED\cdot L}.
\end{align}
\end{subequations}
To further expand on the 1D Conv operation, each step will be highlighted.
First, the input is transposed $U^{\top}_S\in \mathbb{R}^{ED\times L}$ such that the 1D Conv is applied along the sequence length. Furthermore, the input is padded to maintain the input and output size. The padding is added to the input matrix $\tilde{U}$ such that
\begin{equation}
    \tilde{U} = \left[\bar{0}_{l}, U^{\top}_S, \bar{0}_r\right] \in \mathbb{R}^{ED\times L+K-1},
\end{equation}
where for kernel size $K\in \mathbb{N}$ if kernel size $K-1$ is odd, then padding is applied asymmetrically such that $\bar{0}_l\in \mathbb{R}^{ED\times \lfloor(K-1)/2\rfloor}$ and $\bar{0}_r\in \mathbb{R}^{ED\times \lceil\left(K-1\right)/2\rceil}$.
Afterwards, the input $\tilde{U}$ is flattened to a vector 
\begin{equation}
    \bar{u} = \text{vec}\left(\tilde{U}\right)\in \mathbb{R}^{ED(L+K-1)}.
\end{equation} 
Next, because a depthwise 1D Conv is utilized, consider a learnable set of weighted sliding windows $ \kappa \in \mathbb{R}^{ED\times K}$, where for every expanded model dimension $[1, \dots, ED]$ there is a unique sliding window 
\begin{equation}
    \kappa_{d,:} = \left[k_{d,1}, \dots,k_{d,K}\right],
\end{equation}
that is reconstructed into a block diagonal matrix consisting of unique Toeplitz matrices, for every model dimension,  
\begin{subequations}
   \begin{align}
    \mathcal{T} =& \text{Toeplitz}\left( \kappa_{d,:}\right), \, d = 1,\dots, ED, \\
    W_{Conv} =& \text{blk}\left(\mathcal{T}\right).
    \label{eq:ConvToeplitz}   
\end{align} 
\end{subequations}
Lastly, a bias term is added to each model dimension by constructing a vector by repeating each bias $b$, $L$ times, such that
\begin{equation}
    \tilde{b} = \text{col}\left(b, \dots, b\right).
\end{equation}
In the following section, the SSM block will be explained in mathematical detail.

\subsection{SSM}
The SSM function is the core innovation of Mamba, which implements a time-varying content-aware recurrent mechanism for sequence modeling.
Let $\mathbb{F}_{SSM}$ describe the SSM function:
\begin{equation}
   Y_S =  \mathbb{F}_{SSM}\left(U_{\sigma}\right) := H_LC+U_{\sigma}\boxdot D, \label{eq:SSM}
\end{equation}
where the input components are constructed by two sets of operations, i.e., first, the state-space model
\begin{equation}
    H_\tau  = \bar{A}_{\tau}\odot H_{\tau-1}+\bar{B}_{\tau} \boxdot U_{\sigma\tau,:},
\end{equation}
and second, the discretization operations
\begin{subequations}  \label{eq:ZOH}
    \begin{align}
           \bar{\mathcal{A}} &= \exp\left(\Delta_{\tau}\boxdot A\right), \label{eq:ZOH_a}\\
    \bar{\mathcal{B}} &=  \Delta_{\tau}\boxdot  U_{\sigma} W_{B}^{\top}, \\
    \Delta_{\tau} &=  F_{NN}\left(U_{\sigma} W_{\Delta}^{\top}\right), \\ 
            C &=  U_{\sigma} W_{C}^{\top}. 
    \end{align}
\end{subequations}
To elaborate on the SSM function, the state-space modeling and discretization steps are further clarified.

\begin{figure}[h]
    \centering
   \includegraphics[width=0.9\linewidth]{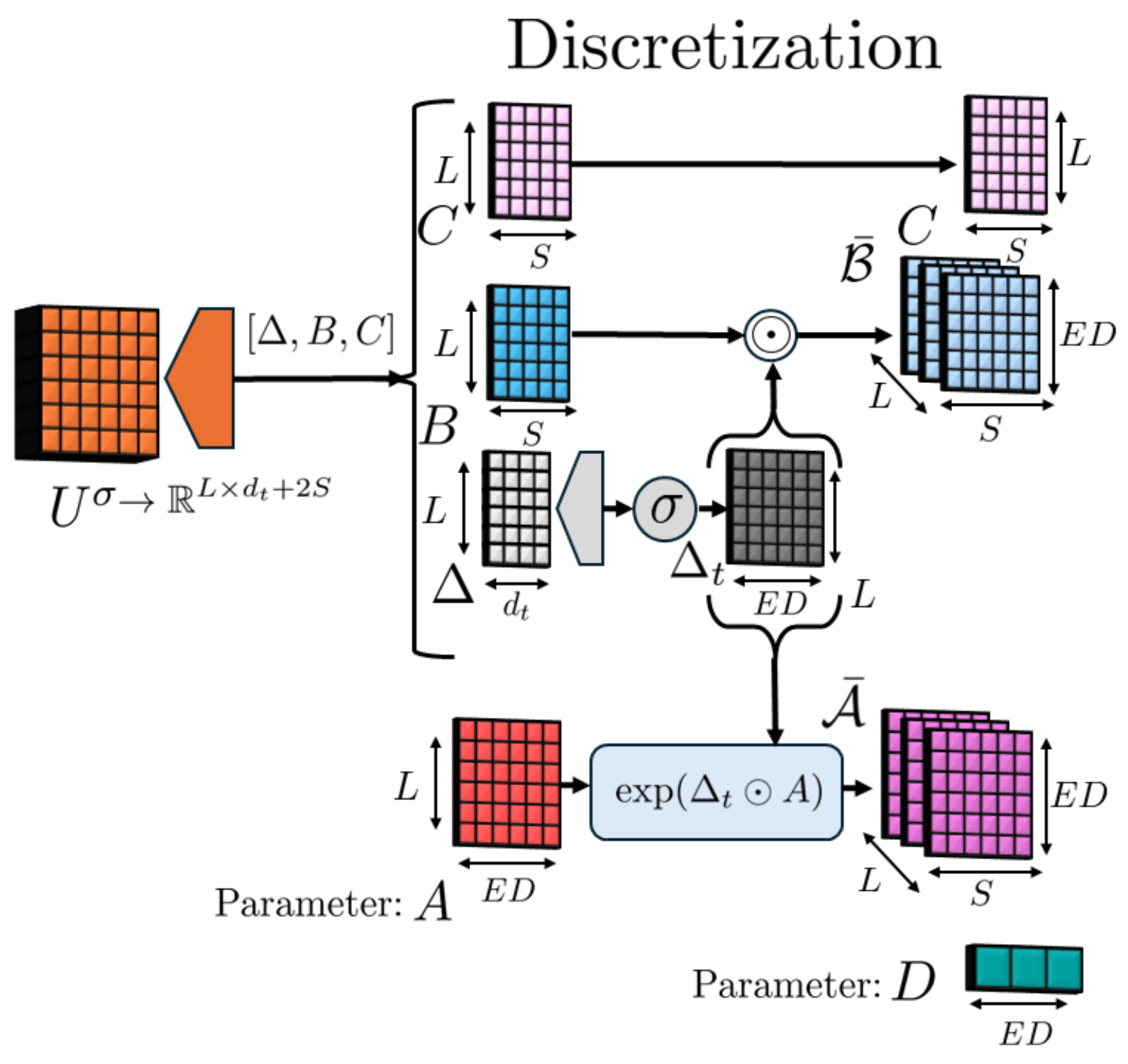}
    \caption{A visual representation of the discretization algorithm in the SSM function.}
    \label{fig:Discretization}
\end{figure}

First, the discretization algorithm is considered,  which incorporates a context aware selection mechanism. This is achieved by constructing $B, C \in \mathbb{R}^{L\times S}$ and $\Delta \in \mathbb{R}^{L \times d_t}$ from the input $U_{\sigma}$ by means of maps $W_{B}^{\top}:  \mathbb{R}^{L \times ED} \to  \mathbb{R}^{L \times S}$, $W_C^{\top}:  \mathbb{R}^{L \times ED} \to  \mathbb{R}^{L \times S} $ and $W_\Delta^{\top}:   \mathbb{R}^{L \times ED} \to  \mathbb{R}^{L \times d_t}$:
\begin{subequations}
    \begin{align}
        B &= U_{\sigma}W_B^{\top}, \label{eq:B}\\
        C &= U_{\sigma}W_C^{\top}, \\
        \Delta &= U_{\sigma}W_ \Delta^{\top}.
    \end{align}
\end{subequations}
By making $B$ and $C$ input dependent, $B$ becomes analogous to an input selection gate and $C$ becomes the output selection gate.
To finalize the construction of the sampling time matrix, a softplus activation is applied, element-wise, to allow $\Delta_t$ to be a nonlinear function that can describe time warping \cite{Tallec2018}.
Furthermore, a map $W_{\Delta_{\tau}}^{\top}: \mathbb{R}^{L \times d_t} \to \mathbb{R}^{L \times ED}$ and bias term $b_{\Delta_{\tau}}$ are added.
Thus, $\Delta_\tau$ is described by:

\begin{equation}
    \Delta_\tau = \ln\left(1+\exp\left(\Delta W_{\Delta_{\tau}}^{\top}+b_{\Delta_{\tau}}\right)\right) =F_{NN}\left(U_{\sigma}\right),
    \label{eq:SamplingTime}
\end{equation}
which can be seen as a feedforward NN.
The application of a softplus activation function is motivated by the connection to the GRU gating mechanism, which is shown in \cite{Gu2023}, where, similarly to $B$ and $C$, by constructing $\Delta_\tau$ from the input, $\Delta_\tau$ becomes analogous to a gate that dynamically focuses on past sequences or current inputs. The temporal focus is achieved by applying discretizations to parameters $A\in \mathbb{R}^{ED\times S}$ and $B$ to construct $\bar{\mathcal{A}}$ and $\bar{\mathcal{B}}$ such that $\bar{\mathcal{A}}$ dictates what information in hidden states $H_\tau$ should be propagated versus what input information should be stored in $H_\tau$ through $\bar{\mathcal{B}}$.
 A graphical illustration is given in Figure \ref{fig:Discretization}.

\begin{figure}
    \centering
    \includegraphics[width=0.8\linewidth]{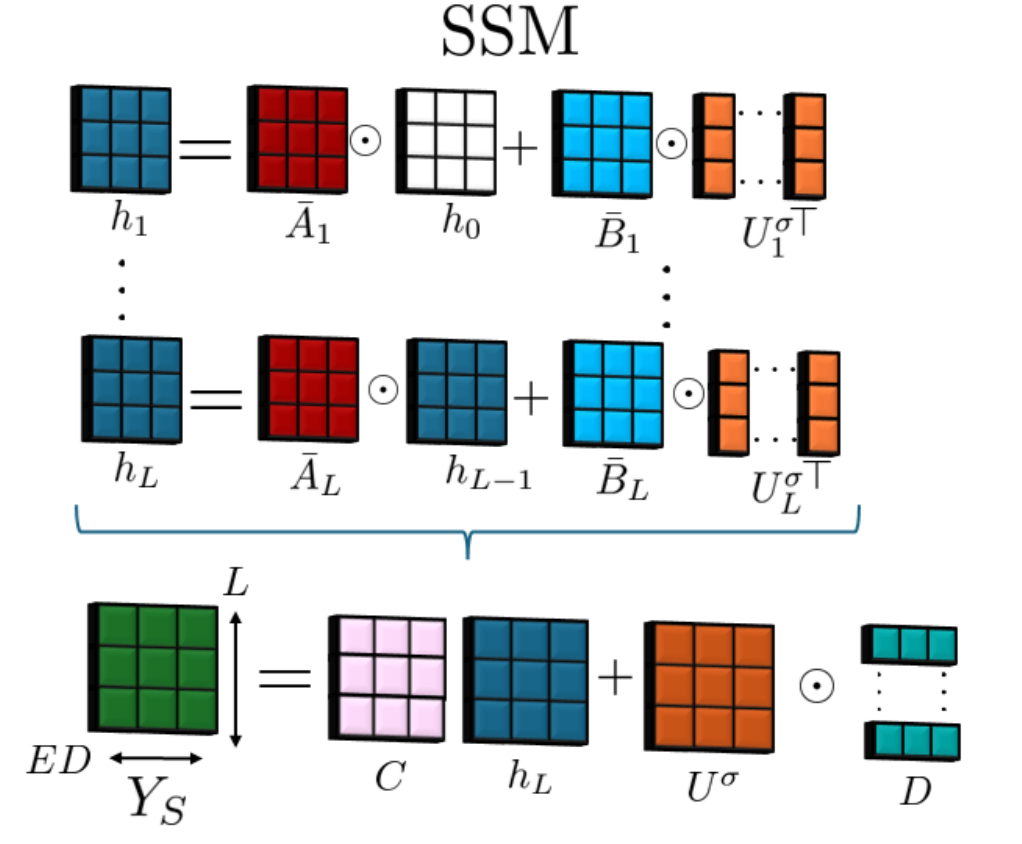}
    \caption{A visual representation of the selective state-space block in a Mamba NN}
    \label{fig:SSM}
\end{figure}

\begin{remark}
It should be noted that $A$ is a dense matrix where the row vectors contain the diagonal elements of individual matrices.
Therefore, tensor $\bar{\mathcal{A}}\left(U_{\sigma}\right)$ is a dense tensor where for every sequence step $l =1, \dots, L$ there is a unique matrix $\bar{A}_l$. Furthermore, every column vector in matrix $\bar{A}_l \in \mathbb{R}^{ED \times S}$ contains the diagonal elements of matrix $\bar{A}_{l,d}$ such that for every sequence step in $L$, there exist $ED$ parallel SSM, as can be seen in Figure \ref{fig:SSM}.
The fact that $A$ is a diagonal matrix is a key component because it allows the ZOH to be approximated using the exponents of the diagonal components rather than a matrix exponential, which would lead to unstable gradients \cite{Lezcano-Casado2019}.
\end{remark}

Next we present the state-space modeling algorithm, which is fundamentally based on \eqref{eq:NLMIMO}. 
For Mamba, a parallelized description of \eqref{eq:NLMIMO} is used which implements the context aware matrices of \eqref{eq:ZOH}. Furthermore, a trainable feed-through parameter $D \in \mathbb{R}^{ED}$ is introduced.  Mamba's parallelized discrete SSM is described by:
\begin{equation}
  H_\tau  = \bar{A}_{\tau}\odot H_{\tau-1}+\bar{B}_{\tau} \boxdot U_{\sigma\tau,:}, \label{eq:MambaDiscretizedSSM}  
\end{equation}
where $\tau = 1,\dots,L$ is fictitious time, $H_0$ is initialized as a zero matrix, and $H_\tau$ can be partitioned as $H_\tau = [h_{\tau,1}, \dots, h_{\tau,ED}] \in \mathbb{R}^{ED\times S}$.
The output of the SSM layer is constructed by:
\begin{equation}
    Y_S = \sum_{s=1}^{S}H_{L_{b,l,d,s}} C_{b,l,s     } + U_{\sigma} \boxdot D. \label{eq:MambaDiscretizedSSM2}
\end{equation}
An illustration of how the SSM mechanism operates is given in Figure \ref{fig:SSM}. This can be regarded as an $L$ step-ahead time-varying linear predictor in the feature space.

\begin{remark}
We emphasize that the parallelized description given in \eqref{eq:MambaDiscretizedSSM} is equivalent to a time-varying linear state-space model.
\end{remark}



\begin{figure*}
    \centering
    \includegraphics[width=1\linewidth]{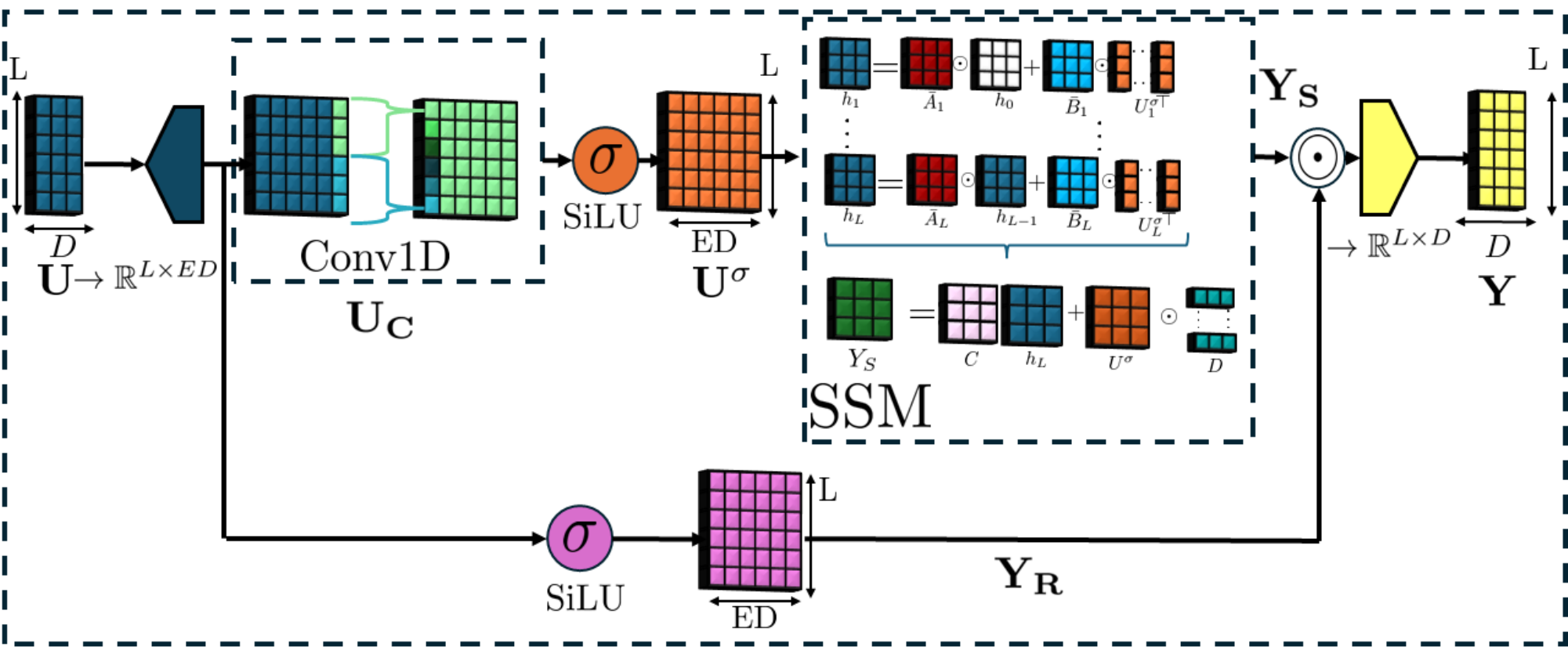}
    \caption{A complete schematic of the Mamba block-wise architecture (for simplicity of illustration, the two different input lifting operators corresponding to the two main parallel paths are merged into one).}
    \label{fig:Mamba_structure}
\end{figure*}

Now that a consistent and complete mathematical description is provided for Mamba, a graphical illustration of the full model is given in Figure \ref{fig:Mamba_structure}. 

In the next section, we will discuss the adaptation required for implementation of Mamba as a multi-step predictor in Model Predictive Control.

\section{Mamba for Predictive Control}\label{sec:Seq2Seq}
In this section, a novel implementation of Mamba for predictive control is introduced that allows Mamba to learn the Seq2Seq mapping for an arbitrary initial state.
To this end, a novel embedding is introduced of $u_{\left [0,N-1 \right]}(k)$ and $x_0(k)$ that maintains the same input sequence length as output sequence length for $y_{\left[1,N \right]}(k)$ such that the usage of a decoder-only architecture is enabled. 

In the previous Section \ref{sec:Mamba}, $L$ was used to describe the length of a sequence, but in control systems it is common to use prediction horizon $N$. From this section onward, the prediction horizon $N$ will be used instead of the sequence length $L$.

\subsection{Data Embedding}
In a decoder-only architecture utilizing Mamba, no changes are made to the sequence length dimension. Therefore, an alternative input embedding can be constructed where the input sequence is concatenated with initial condition columns. 

Consider a dataset $\mathcal{D}=:[X_0,U_f,Y_f]$ with amount of samples $T$
\begin{subequations}
 \begin{align}
    X_0 &= [x_0\left(0\right), \dots, x_0\left(T-1\right)],  \\
U_f &= [u_{\left[0,N-1\right]}\left(0\right)^\top, \ldots, u_{\left[0,N-1\right]}\left(T-1\right)^\top],  \\
Y_f &= [y_{\left[1,N\right]}(1)^\top, \ldots, y_{\left[1,N\right]}\left(T\right)^\top],
\end{align}   
\end{subequations}
where $X_0 \in \mathbb{R}^{T\times 1\times n_x}$, $U_f \in \mathbb{R}^{T\times N\times n_u}$ and $Y_f \in \mathbb{R}^{T\times N\times n_y}$. Then the following input embedding can be constructed for each sample in the dataset.
\begin{definition}{\textit{Input and output embedding for MPC}}
\label{prob:StateEmbedding} 
\begin{equation}
\label{eq:StateEmbed}
\begin{split}
& U_{\mathrm{IC}}\left(k\right)= \mathbb{F}_{\mathrm{IC}}(u_{\left[0,N-1\right]}(k),x_0(k))  :=
\\
&\begin{bmatrix}
u\left(0|k\right)^{\top} & x_0\left(k\right)^{\top} \\
u\left(1|k\right)^{\top} & x_0\left(k\right)^{\top}\\
\vdots &  \vdots \\
u\left(N-1|k\right)^{\top}&x_0\left(k\right)^{\top}
\end{bmatrix} \in \mathbb{R}^{ N\times n_u+n_x}
\end{split}
\end{equation}
where $U_{\mathrm{IC}}\left(k\right)$ is constructed for each sample in dataset $\mathcal{D}$. 
Furthermore, the output embedding for every sample in the batch is given by;
\begin{equation}
Y\left(k\right) := 
    \begin{bmatrix}
        y\left(1|k\right)^{\top}   \\
        \vdots     \\
         y\left(N|k\right)^{\top} 
    \end{bmatrix} \in \mathbb{R}^{ N\times n_y}.
    \label{eq:Yemb}
\end{equation}
\end{definition}
The structure given in \eqref{eq:StateEmbed} allows the embedding of an arbitrary number of initial states, inputs, and outputs and maintains the sequence length.

\subsection{Mamba-MPC}
To effectively construct a Seq2Seq mapping from embedding \eqref{eq:StateEmbed} while utilizing Mamba, adjustments are incorporated.  Let $\mathbb{F}_{\mathrm{MPC}}(U_{\mathrm{IC}})$ denote the function describing the Mamba-MPC function, where $\mathbb{F}_{\mathrm{MPC}}(U_{\mathrm{IC}})$ is described by: 
\begin{equation}
   \hat{Y} =  \mathbb{F}_{\mathrm{MPC}}(U_{\mathrm{IC}}) := \mathbb{F}_{\mathrm{\hat{Y}}} \circ\mathbb{F}_{\mathrm{RMS}}(\mathbb{F}_{\mathrm{Mamba}}\left(U) +U_{\mathrm{in}}, W_{\mathrm{RMS_2}} \right) \label{eq:MambaMPC}
\end{equation}
where
\begin{subequations}
\begin{align}
    U_{\mathrm{in}} &=  \mathbb{F}_{\mathrm{E}}\left(U_{\mathrm{IC}}\right)\\
    U  &= \mathbb{F}_{\mathrm{RMS}}(U_\mathrm{in}, W_{\mathrm{RMS}_1})  \label{eq:InputMPC}.
    \end{align}
\end{subequations}
Although \eqref{eq:MambaMPC} conceptually contains a single $\mathbb{F}_{\mathrm{Mamba}}$, multiple Mamba layers can be cascaded when required. As introduced in \eqref{eq:Fmamba}, the input to Mamba must be of dimension $D$. Therefore, the input embedding is projected to the input dimension $D$ by the map $W_{\mathrm{E}}: \mathbb{R}^{B\times N \times n_u+n_x}\to \mathbb{R}^{B\times N\times D}$ such that:
\begin{equation}
U_{\mathrm{in}} = \mathbb{F}_{\mathrm{E}}(U_{\mathrm{IC}}) := U_{\mathrm{IC}} W_{\mathrm{E}}^{\top} + b_{\mathrm{E}},
\end{equation}
where $b_{\mathrm{E}}$ is a bias term. 

Furthermore, the input to the Mamba block needs to be normalized to allow deeper networks and gradient conditioning \cite{Lubana2021}. The choice was made for 	Root-Mean-Square-Layer-Normalization (RMS norm) \cite{Zhang2019} due to the lower computational cost compared to layer normalization. In \eqref{eq:MambaMPC} there are two RMS norm layers present, which, conceptually, are the same. However, each layer has unique parameters. Generally, the RMS function can be described as:
\begin{equation}
\mathbb{F}_{\mathrm{RMS}}(U,W_{\mathrm{RMS}} ) :=W_{\mathrm{RMS}} \odot \frac{U}{\sqrt{\frac{1}{D} \sum_{k=1}^{D} (U_{:,k})^2 + \epsilon}}
\end{equation}
where $W_{\mathrm{RMS}}$ is the weighting, which is a learnable parameter, and $\epsilon$ is a small static bias.
Furthermore, a residual connection from the output of $\mathbb{F}_{\mathrm{E}}$ is added to improve gradient flow.
Therefore, the output of the second RMS norm is given by:
\begin{equation}
    Y_{\mathrm{RMS}} = \mathbb{F}_{\mathrm{RMS}}\left(Y_{\mathrm{Mamba}}+U_\mathrm{in}, W_{\mathrm{RMS_2}}\right)
\end{equation}

\begin{remark}
Although residual connections are commonly used for gradient flow, residual connections also add linear time invariant (LTI) components to the output which simplifies the learning problem. Therefore, Mamba-MPC can be separated into an LTI component, a nonlinear time invariant component in the Mamba residual branch \eqref{eq:Residual}, and a temporal modeling branch in the SSM block \eqref{eq:SSM}. 
\end{remark}

Lastly, a linear projection is utilized to map the model dimension to the output dimension by means of map $W_{\hat{Y}}: \mathbb{R}^{B \times N \times D} \to \mathbb{R}^{B \times N \times n_y}$ :
\begin{equation}
    \hat{Y} =  Y_{\mathrm{RMS}}W_{\hat{Y}}^{\top}+b_{\hat{Y}},
\end{equation}
where $b_{\hat{Y}}$ is a bias term.

Now that the network structure of Mamba-MPC has been defined,  a model learning problem can be constructed:

\begin{problem}{\textit{Mamba-MPC learning problem}}
\label{prob:MambaPredictionModel}
\begin{subequations}\label{eq:LearningProblem}
\begin{align}
&\arg  \min  J\left(Y^{\text{train }}, \hat{Y}^{\text{train }}\right)  = \frac{||Y^{\text {train }}-\hat{Y}^{\text {train }}||_2^2 }{||Y^{\text {train }}||_2^2} \label{eq:RSE}\\
&\text { subject to: } \nonumber\\
&\hat{Y}^{\text{train}}= \mathbb{F}_{\mathrm{MPC}} \circ \mathbb{F}_{\mathrm{IC}}\left(U_f^{\text {train}},X_{0}^{\text {train}}\right) .\label{eq:PredictionModel}
\end{align}
\end{subequations}
\end{problem}

Above, $U_f^{\text {train}}$ is the set of input sequences, $X_{0}^{\text {train}}$ is the set of initial conditions and $Y^{\text{train }}$ is the set output sequences, when used in training.
In conjunction, checkpointing is applied to avoid overfitting and ensure sufficient exploration of the loss landscape. 

Next, the prediction model can be used in the predictive control for tracking a reference $r(k)= [r(0|k), \dots, r(N|k)]$, with the error between the output sequence and the reference defined as $\varepsilon_{y_{[0,N]}}:=y_{[0,N]}(k)-r_{[0,N]}(k)$. Furthermore, the change in input is penalized given by $\Delta u_{[0,N-1]}(k) = u_{[0,N-1]}(k)-u_{[-1,N-2]}(k)$, where $u(-1|k):=u(k-1)$.  
\begin{problem}{\textit{Mamba-MPC}}
\label{prob:Mamba-MPC} 
\[
\begin{split}
&\min_{\begin{matrix}u_{[0,N-1]}(k)\\ y_{[1,N]}(k)\end{matrix}} \varepsilon_{y_{[0,N]}}^{\top} \Omega\varepsilon_{y_{[0,N]}} 
 + \Delta u_{[0,N-1]}(k)^{\top} \Psi \Delta u_{[0,N-1]}(k) \\
& \text{ subject to: } \\
&y_{[1,N]}(k) =  \mathbb{F}_{\mathrm{MPC}} \circ \mathbb{F}_{\mathrm{IC}}\left(u_{\left[0,N-1\right]}(k),x_0\left(k\right)\right), \\
&u_{\left[0,N-1\right]}\left(k\right)\in\mathbb{U}\times\ldots\times\mathbb{U},\\
&y_{\left[1,N\right]}\left(k\right)\in\mathbb{Y}\times\ldots\times\mathbb{Y},\\
&y(0|k)=y(k).
\end{split}
\]
\end{problem}
Above the quadratic cost function weight matrices are defined as follows:
$$
\Omega=\left[\begin{array}{llll}
Q & & & \\
& \ddots & & \\
& & Q & \\
& & & P
\end{array}\right], \quad \Psi=\left[\begin{array}{lll}
R & & \\
& \ddots & \\
& & R
\end{array}\right].
$$ 
Given that $\mathbb{F}_{\mathrm{MPC}}$ is constructed by Mamba, this problem is considered non-convex while being locally Lipschitz and globally differentiable. Therefore, non-linear optimizers such as IPOPT and SQP solvers will be considered.
In the next section, we analyze the performance of Mamba-MPC as defined in Problem~\ref{prob:Mamba-MPC} by means of 2 nonlinear benchmark examples.

\subsection{Towards Stability Guarantees for Mamba-MPC}
Proving stability of Mamba-MPC is non--trivial due to the complex heterogeneous Mamba architecture, which includes convolution operations, time-varying discretization, recurrent linear models and Hadamard products of feedforward neural network outputs. A promising framework for addressing stability of MPC with general input--output predictors, such as neural networks, was recently developed in \cite{Bold2026Lipschitzstability}. Therein, exponential stability is established without terminal ingredients under the assumptions that both the prediction model and the system are Lipschitz continuous, and that the predictor satisfies proportional error bounds. These results could be applied to Mamba--MPC provided that Lipschitz bounds can be established for the Mamba NN architecture.

While a complete stability analysis is beyond the scope of this paper, in what follows we point out how Lipschitz bounds can be derived for the Mamba network. To this end we first introduce the definition of Lipschitz continuity tailored to our needs.

\begin{definition}
A map $\mathbb{F}:\mathbb{R}^n\times \mathbb{R}^m\rightarrow\mathbb{R}^p$ is Lipschitz continuous
w.r.t. its first argument uniformly in its second argument if there exists a Lipschitz constant $L>0$, such that
$$
\|\mathbb{F}(x,u) - \mathbb{F}(x',u)\| \leq L\|x-x' \|, \quad \forall u\in\mathbb{U}
$$
for all $x,x'\in\mathbb{X}$.
\end{definition}
Next, establishing Lipschitz continuity of the Mamba network amounts to verifying that each of its underlying functions/maps is Lipschitz continuous. In particular, the network employs the SiLU activation function \eqref{eq:sigma} and the Softplus activation function \eqref{eq:SamplingTime}, both of which are known to be globally Lipschitz (see, e.g., \cite{sun2025entropy, luo2026principles}). Consequently, these nonlinearities do not obstruct the Lipschitz continuity of the overall network.

Another operation performed in the network is the discretization of the state-space equation \eqref{eq:ZOH_a}, implemented in the SSM block and involving a matrix exponential term. This exponential is applied element-wise. It is well known that exponential functions are not globally Lipschitz on \(\mathbb{R}\). However, they are locally Lipschitz, which can be shown directly. Let \(f(x)=e^x\), and consider any bounded interval \([a,b]\subset \mathbb{R}\). For any \(x,y\in [a,b]\), the Mean Value Theorem guarantees the existence of some \(c\) between \(x\) and \(y\) such that $e^x - e^y = e^c (x - y)$.
Taking absolute values yields
\[
|e^x - e^y| = e^c |x - y| \le e^b |x - y|
\]
where $b\geq c$. Thus, \(e^x\) is Lipschitz continuous on every bounded interval, with Lipschitz constant \(e^b\). Consequently, \(e^x\) is locally Lipschitz on \(\mathbb{R}\).

Another concern that might break the Lipschitz continuity of the Mamba network is the case when the state-space model \eqref{eq:MambaDiscretizedSSM} is unstable, and states might grow unbounded. Since the state space model is a time varying system, stability analysis is not straightforward. This time varying behavior is introduces by the input dependent discretization \eqref{eq:ZOH}. However, recently it was shown in \cite{cao2026a} that the SSM recursion \eqref{eq:MambaDiscretizedSSM} is BIBO stable under the mild assumption that the matrix $A$ in \eqref{eq:ZOH_a} is Hurwitz.

Finally, the outputs of the SSM block \(Y_S\) and the parallel path \(Y_R\) are combined via element-wise multiplication, i.e., the Hadamard product. Even if both mappings are individually Lipschitz, this alone does not guarantee that their Hadamard product is Lipschitz. Additional assumptions, such as boundedness of the components, are required to ensure Lipschitz continuity of the combined map.
\begin{lemma}
Let \(f,g:\mathbb{R}\to\mathbb{R}\) be Lipschitz and bounded, i.e., there exist constants \(L_f, L_g >0\) such that:
\[
|f(x)-f(y)| \le L_f |x-y|, \quad |g(x)-g(y)| \le L_g |x-y|.
\]
Then \(h(x)=f(x)g(x)\) is also Lipschitz. 
\end{lemma}
\begin{proof} Let \(h(x)=f(x)g(x)\). Then
$$
h(x)-h(y)=f(x)(g(x)-g(y)) + g(y)(f(x)-f(y)).
$$
Taking absolute values and applying the triangle inequality gives
$$
|h(x)-h(y)| \le |f(x)|\,|g(x)-g(y)| + |g(y)|\,|f(x)-f(y)|.
$$
Using the Lipschitz bounds and boundedness of $f(x)\leq M_f$ and $g(x)\leq M_g$ for all $x$, we obtain
$$
|h(x)-h(y)| \le (M_f L_g + M_g L_f)\,|x-y|.
$$
Hence \(h\) is Lipschitz.
\end{proof}

The above analysis provides the key steps for establishing local Lipschitz boundedness of the complete Mamba network architecture. This provides a useful basis for future work that will deal with a rigorous stability analysis of Mamba--MPC, for instance within the framework of \cite{Bold2026Lipschitzstability}.


\section{Validation in Numerical Examples}\label{sec:NumExample}
This section details the comprehensive evaluation of the proposed Mamba-MPC scheme across a series of performance metrics.

For all considered systems, we evaluate the performance of Mamba-MPC against a state-space LSTM \cite{awryczuk2025} based MPC controller as in \cite{deJong2025}. For the Van der Pol oscillator an additional comparison with a Transformer-MPC model was made. The Transformer-MPC model was constructed according to \cite{Park2023} with the addition of a temperature parameter in the self-attention mechanism. The addition of a temperature parameter was done to aid numerical stability \cite{emadi2026}. The utilized parameters can be found in Table~\ref{tab:VDP_Results} and Table \ref{tab:FourTank_Results} for the Van der Pol and Four tank, respectively.

\begin{table}[h!]
\centering
\caption{Comparison of Mamba-MPC and LSTM-MPC for the Van der Pol System.}
\label{tab:VDP_Results}
\resizebox{\columnwidth}{!}{
\begin{tabular}{l cccc}
\toprule
& \multicolumn{4}{c}{\textbf{Van der Pol}} \\
\cmidrule(lr){2-5}
& \multicolumn{2}{c}{Reference Tracking} & \multicolumn{2}{c}{Noise} \\
\cmidrule(lr){2-3} \cmidrule(lr){4-5}
& Mamba & LSTM & Mamba & LSTM \\
\midrule
\multicolumn{5}{l}{\textit{Model Architecture}} \\
Model Dim. $D$            & 8 & 2 & 8 & 2 \\
Kernel Size               & 10 & - & 10 & - \\
State Dim.                & 8 & 26 & 8 & 26 \\
Horizon $N$               & 10 & 10 & 10 & 10 \\
Layers                    & 6 & - & 6 & - \\
Params                    & 3418 & 3052 & 3418 & 3052 \\
\midrule
\multicolumn{5}{l}{\textit{Controller Performance}} \\
MAE                       & $\color{mygreen} 0.066$ & $0.072$ & $\color{mygreen} 0.07$ & $0.09$ \\
MSE                       & $\color{mygreen} 0.058$ & $0.066$ & $\color{mygreen} 0.03$ & $0.04$ \\
Comp. Time (s)            & $\color{mygreen} 0.018$ & $0.20$ & $\color{mygreen} 0.07$ & $0.71$ \\
Sampling Time (s)         & 0.1 & 0.1 & 0.1 & 0.1 \\
\midrule
\multicolumn{5}{l}{\textit{Training Metrics}} \\
Time (min)                & 400 & 35 & 400 & 35 \\
Train Loss                & $2 \cdot 10^{-5}$ & $8.4\cdot 10^{-7}$ & $6.8 \cdot 10^{-5}$ & $6.2 \cdot 10^{-4}$ \\
Val. Loss                 & $5.5 \cdot 10^{-5}$ & $6.9\cdot 10^{-5}$ & $8.8\cdot 10^{-5}$ & $2.2 \cdot 10^{-2}$ \\
\bottomrule
\end{tabular}%
}
\end{table}

\begin{table}[h!]
\centering
\caption{Transformer-MPC.}
\label{tab:Transformer_Results}
{
\begin{tabular}{l c}
\toprule
& \multicolumn{1}{c}{\textbf{Van der Pol}} \\
\cmidrule(lr){2-2}
& \multicolumn{1}{c}{Reference Tracking}  \\
\cmidrule(lr){2-2} \cmidrule(lr){2-2}
& Transformer \\
\midrule
\multicolumn{2}{l}{\textit{Model Architecture}} \\
Model Dim. $D$            & 8  \\
Heads               & 4 \\
MLP Dim.                & 32 \\
Horizon $N$               & 10  \\
Layers                    & 3  \\
Params                    & 26600  \\        
\midrule
\multicolumn{2}{l}{\textit{Controller Performance}} \\
MAE                       & $0.13$  \\
MSE                       & $0.15$  \\
Comp. Time (s)            & $0.09$  \\
Sampling Time (s)         & 0.1  \\
\midrule
\multicolumn{2}{l}{\textit{Training Metrics}} \\
Time (min)                & 1200  \\
Train Loss                & $7 \cdot 10^{-5}$ \\
Val. Loss                 & $2.3\cdot 10^{-4}$ \\
\bottomrule
\end{tabular}%
}
\end{table}

First, the SISO Van der Pol oscillator is considered where stabilization is tested for multiple initial conditions. Subsequently, a tracking performance comparison against LSTM-MPC and Transformer-MPC is made where a fair assessment is ensured by using identical datasets, MPC cost function, prediction horizon and each network is optimized for their respective optimal compromise between accuracy, size and numerical accuracy. Thereafter, a comparison is made for computational scaling against prediction horizon. Then, 
robustness against noise is tested during reference tracking and compared against LSTM-MPC. Afterwards, the MIMO Four Tank benchmark is introduced where reference tracking for multiple outputs is tested and compared with NMPC, where again, a fair comparison is made in terms of MPC cost and prediction horizon. Finally, the Mamba-MPC is applied to a physical benchmark system, the Quanser Aero2, to validate real-time reference tracking performance.
 All results presented in this section were generated using PyTorch and CasADI \cite{Andersson2019} frameworks which were executed on an Nvidia RTX 4070 and AMD Ryzen 9 7845HX CPU, respectively.

\subsection{Van der Pol Oscillator}
\begin{figure}
    \centering
    \includegraphics[width=1\linewidth]{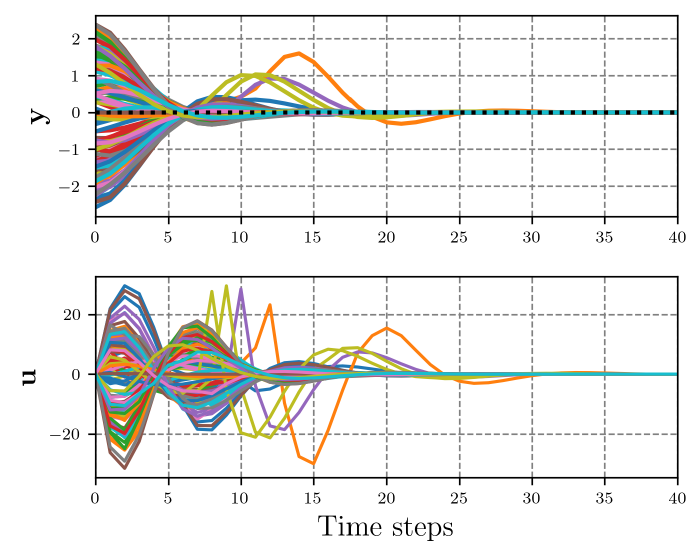}
    \caption{Closed-loop performance of Mamba-MPC for stabilization of a Van der Pol oscillator for $100$ different initial conditions.}
    \label{fig:InitialConditions}
\end{figure}

Mamba-MPC is first evaluated on the SISO Van der Pol system with control input, i.e., 
\begin{subequations}
    \begin{align}
\dot{x}_1 & =x_2, \\
\dot{x}_2 & =\mu\left(1-x_1^2\right) x_2+u, \\
y & =x_1,
\end{align}
\end{subequations}
with $\mu=1$, input $u\in \mathbb{R}$, output $y\in \mathbb{R}$ and initial states $x_0\in \mathbb{R}^{2}$. Furthermore, the Van der Pol is discretized with forward Euler with sampling time $T_s=0.1$ s. The specific selection of $\mu=1$ ensures that Van der Pol exhibits limit cycle behavior, thereby constituting a non-trivial benchmark problem. A dataset consisting of $T=40000$ samples is created via an open-loop identification experiment, in which a multisine input signal is used featuring a peak amplitude of $15$, across a frequency range of approximately $0.0049$ Hz to $4.88$ Hz, containing $30$ harmonics.
A training and validation dataset is created according to \eqref{eq:StateEmbed} and utilized in learning Problem \eqref{eq:LearningProblem}. A model is constructed according to the model parameters given in Table \ref{tab:VDP_Results}. 
The model is trained for $4000$ epochs in PyTorch, using the Adam optimizer \cite{Kingma2014} with an initial learning rate $l_r = 1\cdot 10^{-3}$, an L2 regularization cost of $\lambda = 10^{-5}$, weight decay factor of $w_t=1\cdot10^{-5}$, and exponential learning rate scheduler with $\gamma=0.998$ per $10$ epochs. 
\\
\\
The PyTorch model is implemented as a CasADI function to formulate the prediction model in \ref{prob:Mamba-MPC} to enable closed-loop control. 
We first assess the system's ability to stabilize a system for a variety of initial conditions. To generate a diverse set of initial conditions, $x_0 =\text{col}(x_1(0), x_2(0))$, $100$ realizations are drawn according to $x_{1}(0)\sim\text{Uniform}(-2.5, 2.5)$ and $x_2(0)\sim\text{Uniform}(-2.0, 2.0)$. For the Mamba-MPC scheme \eqref{prob:Mamba-MPC} the prediction horizon is $N=10$, $Q=50$, $R=0.5$ and $P=100$. This optimization problem is solved using IPOPT solver \cite{Wachter2006}, which is employed for all Van der Pol examples. Consider the stabilizing closed-loop response of Mamba-MPC for the Van der Pol for $100$ initial condition realizations in Figure \ref{fig:InitialConditions}. The results confirm that Mamba-MPC successfully stabilizes the system for all sampled trajectories, though one specific trajectory exhibited a higher amplitude of oscillations during the convergence phase.

To evaluate the tracking performance of Mamba-MPC, a comparative analysis is conducted against a Transformer and LSTM model. Mamba-MPC, Transformer, and LSTM are trained on a identical dataset and are trained according to learning Problem \eqref{eq:LearningProblem}. The choice of the LSTM as baselines is motivated by the architectural similarities between Mamba and LSTM specifically regarding their recurrent states and gating mechanism. The Transformer is chosen due to having comparable to superior modeling performance in comparison to Mamba in a wide range of modalities. As previously mentioned, a temperature parameter was added to the self-attention of the Transformer and was empirically chosen to be $2.5$.
The Transformer settings can be found in Table \ref{tab:Transformer_Results}. To ensure an equitable comparison, the hyperparameter setting is chosen to be optimal for its architecture in terms of online computational time, modeling performance, and online numerical stability. All control schemes are implemented with a prediction horizon $N=10$ and $Q=100$, $R=0.5$.

\begin{figure}
    \centering
    \includegraphics[width=1\linewidth]{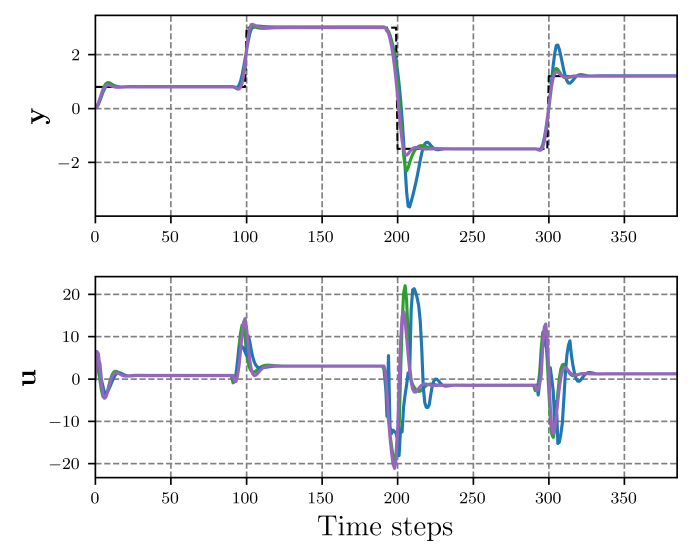}
    \caption{Comparison between Mamba-MPC (\textcolor{mypurple}{\rule[0.5ex]{1em}{1pt}}), LSTM-MPC (\textcolor{mygreen}{\rule[0.5ex]{1em}{1pt}}), and Transformer-MPC (\textcolor{myblue}{\rule[0.5ex]{1em}{1pt}})  for a piece-wise constant reference (\textcolor{myblack}{\rule[0.5ex]{0.2em}{1pt}} \textcolor{myblack}{\rule[0.5ex]{0.2em}{1pt}}) on a Van der Pol oscillator.}
    \label{fig:MambaVsLSTM}
\end{figure}

\begin{figure}
    \centering
    \includegraphics[width=0.8\linewidth]{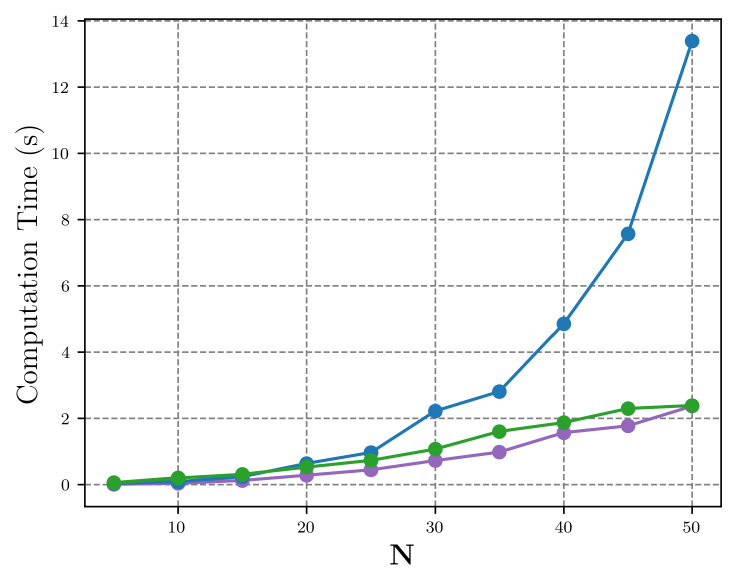}
    \caption{Comparison between the mean online computation time of Mamba-MPC (\textcolor{mypurple}{\rule[0.5ex]{1em}{1pt}}), LSTM-MPC (\textcolor{mygreen}{\rule[0.5ex]{1em}{1pt}}), and Transformer-MPC (\textcolor{myblue}{\rule[0.5ex]{1em}{1pt}}) on a Van der Pol oscillator.}
    \label{fig:CompTime}
\end{figure}

Consider the closed-loop response for tracking a piece-wise constant reference signal in Figure \ref{fig:MambaVsLSTM}. Generally, the LSTM and Mamba based MPC controllers demonstrate a comparable ability to track the reference signal, a finding substantiated when considering the Mean-Absolute-Error (MAE) and Mean-Squared-Error (MSE) in Table \ref{tab:VDP_Results} and table \ref{tab:Transformer_Results}. However, the Transformer-MPC shows significantly worse performance, particularly in regions with aggressive reference changes. In general to obtain a result, the Transformer had to be trained for $N=50$ to generalize beyond a reference level of $2.5$.
We hypothesize that this is due to the numerical difficulties as well as poor generalization and requires a dedicated analysis to form a conclusion. 
However, when considering the mean computation time, Mamba-MPC outperforms both competitors with a computation time of $0.018$s, whereas LSTM-MPC does not stay within the sampling time of the system at $0.2$s, and the Transformer marginally stays within the sampling time at $0.09$s. To further analyze the computational scaling of each architecture, a comparison was made for a set of prediction horizons $N=[5,10,\dots, 50]$. It is evident from Figure \ref{fig:CompTime} that the computational scaling of each model is clearly visible where Mamba and LSTM scale linearly with prediction horizon while the Transformer scales quadratically. Therefore, for following examples, we omit the Transformer since it would be an infeasible approach. 

\begin{figure}
    \centering    
    \includegraphics[width=1\linewidth]{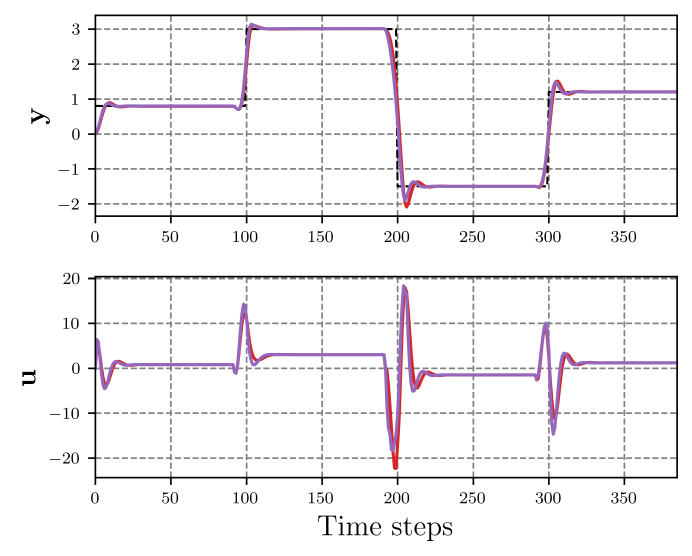}
    \caption{Comparison between stateless Mamba-MPC (\textcolor{mypurple}{\rule[0.5ex]{1em}{1pt}}) and stateful Mamba-MPC (\textcolor{myred}{\rule[0.5ex]{1em}{1pt}}) for a piece-wise constant reference on the Van der Pol in noiseless conditions.}
    \label{fig:StatefulvsStateless}
\end{figure}

\begin{remark}
In section \ref{sec:Seq2Seq}, a choice was made to reset the hidden state at each time step. This was done such that a consistent comparison of Mamba as input-output operator could be constructed. However, the question whether this is optimal remains open. To construct a brief empirical analysis between stateful and stateless training, a comparison was conducted. Identical settings as provided in Table \ref{tab:VDP_Results} we utilized for the stateless and stateful model. To train the stateful model, an adaptation to the training and test set were done to align subsequent sequences along the batch dimension. Thus, the model can propagate $h(0)$ for each batch and maintain chronological order. The training was done for a noiseless data set. In Figure \ref{fig:StatefulvsStateless} the comparison can be seen. It can be noted that the performance is comparable with small deviations. However, to form a strong conclusion a broader study is required. 
\end{remark}

 Lastly, we evaluate the reference tracking capabilities of Mamba-MPC under stochastic noise conditions. We perform a direct comparison with LSTM-MPC. The dataset acquired in the previous example is corrupted with a Gaussian, zero-mean, additive white noise realization such that the measured states $x_0$ and output $y$ have a signal-to-noise ratio of $SNR=20$. Subsequently, both the Mamba and LSTM models are retrained according to \eqref{eq:LearningProblem} with the same settings as in the previous example.

\begin{figure}
    \centering
    \includegraphics[width=1\linewidth]{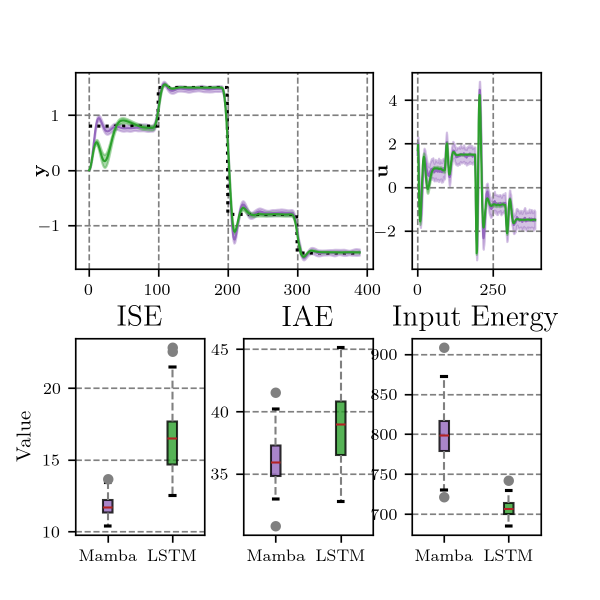}
    \caption{A closed-loop reference (\textcolor{myblack}{\rule[0.5ex]{0.2em}{1pt}} \textcolor{myblack}{\rule[0.5ex]{0.2em}{1pt}}) tracking  comparison between Mamba-MPC  (\textcolor{mypurple}{\rule[0.5ex]{1em}{1pt}}) and LSTM-MPC (\textcolor{mygreen}{\rule[0.5ex]{1em}{1pt}}) under noise conditions with a  Signal-to-Noise Ratio (SNR) conditions of 20. The solid line indicates the mean value, and the shaded area represents the standard deviation obtained across 100 noise realizations.}
    \label{fig:SNR20}
\end{figure}
 
 In closed-loop control implementation, zero-mean Gaussian white noise with standard deviation $s = [0.16, 0.13]$ is added to each measured state, respectively. A total of $100$ independent noise realizations are created to perform a Monte-Carlo evaluation. Both control schemes are implemented with a prediction horizon $N=10$ and quadratic cost $Q=50$, $R=1$ and terminal cost $P=10$. We quantitatively compare the control schemes for the following criteria: 
 
\begin{align}
&\operatorname{ISE}  =\sum_{k= 0 }^{t_{\text {max }}}\|y(k)-r(k)\|_2^2, \quad \text { IAE }=\sum_{k=0}^{t_{\text {max }}}\|y(k)-r(k)\|_1, \nonumber\\
&\text { Input Energy }  =\sum_{k= 0}^{t_{\text {max }}}\|u(k)\|_2^2,
\end{align}

where $\operatorname{ISE}$ indicates the integral squared error and $\text {IAE}$ indicates the integral absolute error.
The resulting mean tracking performance and associated standard deviations are presented in Figure \ref{fig:SNR20} indicating that Mamba-MPC effectively mitigates measurement noise effects. Quantitatively, the $\text{ISE}$ and $\text{IAE}$ suggest that Mamba-MPC achieves superior tracking performance compared to LSTM-MPC, albeit demanding a higher magnitude of input energy. However, closer inspection of the tracking performance, reveals that Mamba-MPC generally exhibits a higher standard deviation than LSTM-MPC. Furthermore, LSTM-MPC maintains a lower tracking error across the simulation, excluding the start of simulation. This initial suboptimality in LSTM-MPC is likely due to its reliance on the hidden state for predictions, where the initialization of the hidden state at initial time $t=0$ is non-trivial. Once the hidden state converges, LSTM-MPC demonstrates a higher degree of noise robustness, in comparison to Mamba-MPC. We suspect that this is due to reinitialization of the hidden states of Mamba-MPC as $0$ at every control interval. Consequently, this renders Mamba-MPC more sensitive to perturbations in the initial condition. Nevertheless, when considering the mean computation time in \ref{tab:VDP_Results}, Mamba-MPC is significantly faster than LSTM-MPC with a mean computation time of $0.07$s and $0.71$s, respectively. This indicates that in the presence of noise, Mamba-MPC is a reliable candidate.

\subsection{Four Tank System}
As a second example, a system with four interconnected tanks (Four Tank) is considered from \cite{Blaud2022} as a benchmark for MIMO systems, i.e, 
\begin{subequations}
    \begin{align}
& \dot{x}_1=\frac{-a_1}{S_c} \sqrt{2 g x_1}+\frac{a_3}{S_c} \sqrt{2 g x_3}+\frac{\gamma_a}{3600 S_c} u_1, \label{eq:FTx1_dot} \\
& \dot{x}_2=\frac{-a_2}{S_c} \sqrt{2 g x_2}+\frac{a_4}{S_c} \sqrt{2 g x_4}+\frac{\gamma_b}{3600 S_c} u_2, \label{eq:FTx2_dot} \\
& \dot{x}_3=\frac{-a_3}{S_c} \sqrt{2 g x_3}+\frac{\left(1-\gamma_b\right)}{3600 S_c} u_2, \label{eq:FTx3_dot} \\
& \dot{x}_4=\frac{-a_4}{S_c} \sqrt{2 g x_4}+\frac{\left(1-\gamma_a\right)}{3600 S_c} u_1 , \label{eq:FTx4_dot}
\end{align}
\end{subequations}
where each state $col(x_1,x_2,x_3,x_4)$ represents the water level of each of the four tanks, respectively. The parameters of the Four Tank system are given in Table \ref{tab:FourTank}. For the Four Tank example we assume that all states are measurable. Note that the system outputs and initial conditions are the states $x,x_0 \in\mathbb{R}^4$ and the inputs are $u \in \mathbb{R}^{2}$, hence a MIMO system. To gather data, an open-loop experiment is performed using a piece-wise-constant input signal with sampling time $T_s =5$ and $T=80000$ samples of input-state data are collected. The dataset is generated according to \eqref{eq:StateEmbed} and used in training according to \eqref{eq:LearningProblem}. The model is trained for a total of $3000$ epochs with the Adam optimizer \cite{Kingma2014} with identical hyperparameters and settings as in the Van der Pol example.
A model is trained with parameters given in \ref{tab:FourTank_Results}.

\begin{table}
\centering
\caption{Comparison of Mamba-MPC and LSTM-MPC for the MIMO Four Tank System.}
\label{tab:FourTank_Results}
\begin{tabular}{l c c}
\toprule
& \multicolumn{2}{c}{\textbf{Four Tank (MIMO)}} \\
\cmidrule(lr){2-3}
& Mamba & LSTM \\
\midrule
\multicolumn{2}{l}{\textit{Model Architecture}} \\
Model Dim. $D$            & 8 & 6 \\
Kernel Size               & 20 & - \\
State Dim.                & 8 & 20 \\
Horizon $N$               & 20 & 20 \\
Layers                    & 1 & - \\
Params                    & 708 & 2324 \\
\midrule
\multicolumn{2}{l}{\textit{Controller Performance}} \\
MAE                       & $\color{mygreen}[2, 1, 1, 1]\cdot10^{-2}$ & $[4, 4, 5, 5]\cdot10^{-2}$ \\
MSE                       & $\color{mygreen}[4, 3, 1, 1]\cdot10^{-3}$ & $[4, 5, 5, 7]\cdot10^{-3}$ \\
Comp. Time (s)            & $\color{mygreen} 0.14$ & $0.65$ \\
Sampling Time (s)         & 5.0 & 5.0 \\
\midrule
\multicolumn{2}{l}{\textit{Training Metrics}} \\
Time (min)                & 400 & 180 \\
Train Loss                & $5.1 \cdot 10^{-7}$ & $5.1 \cdot 10^{-5}$ \\
Val. Loss                 & $1.1 \cdot 10^{-5}$ & $3.2 \cdot 10^{-4}$ \\
\bottomrule
\end{tabular}
\end{table}

\begin{table}
\centering
\caption{Parameters for the Four-Tank System}
\begin{tabular}{|l|l|}
\hline
\textbf{Parameter} & \textbf{Value} \\
\hline
Tank cross-sectional area ($S_c$) & $0.06 \text{ m}^2$ \\
Leakage orifice area ($a_1$) & $1.31 \times 10^{-4} \text{ m}^2$ \\
Leakage orifice area ($a_2$) & $1.51 \times 10^{-4} \text{ m}^2$ \\
Leakage orifice area ($a_3$) & $9.27 \times 10^{-5} \text{ m}^2$ \\
Leakage orifice area ($a_4$) & $8.82 \times 10^{-5} \text{ m}^2$ \\
Three-way valve opening ratio ($\gamma_a$) & $0.3$ \\
Three-way valve opening ratio ($\gamma_b$) & $0.4$ \\
Gravitational acceleration ($g$) & $9.81 \text{ m/s}^2$ \\
\hline
\end{tabular}
\label{tab:FourTank}
\end{table}

\begin{figure}[h]
    \centering
    \includegraphics[width=1\linewidth]{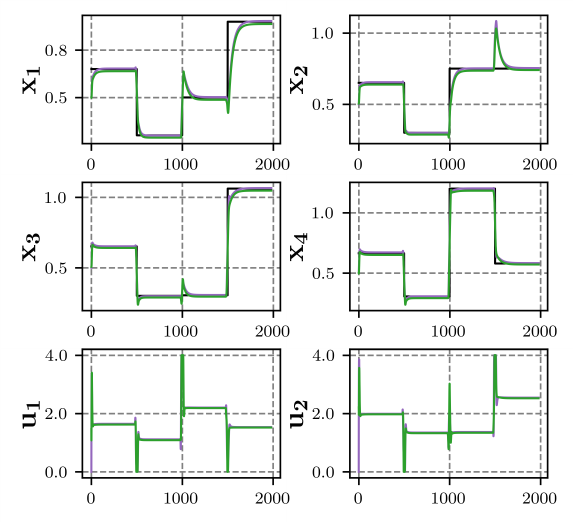}
    \caption{Closed-loop performance comparison between Mamba-MPC (\textcolor{mypurple}{\rule[0.5ex]{1em}{1pt}}) and LSTM-MPC (\textcolor{mygreen}{\rule[0.5ex]{1em}{1pt}}) on the MIMO Four Tank system for piece-wise constant reference trajectories (\textcolor{myblack}{\rule[0.3ex]{1em}{1pt}}).}
    \label{fig:FourTank-MPC}
\end{figure}

Similarly to the Van der Pol example, an evaluation of Mamba-MPC is made against LSTM-MPC. The models are trained on identical datasets according to learning Problem \eqref{eq:LearningProblem} where the parameter budget is ideally set to be approximately equal. However, for LSTM-MPC, to obtain an accurate model a higher amount of parameters were necessary, therefore the parameter budget was ignored.
The trained models are implemented into the MPC scheme \eqref{prob:Mamba-MPC} which is constructed with input constraints $\mathbb{U} = \{u\in R^2| col(0,0)\leq u \leq col(4,4)\}$, tracking cost $Q=100$ and input cost $R = 1$. Consider the control input and closed-loop output given in Figure \ref{fig:FourTank-MPC}, it can be  seen that Mamba-MPC is able to track the piece-wise constant references for all system states and significantly outperforms LSTM-MPC whereas containing a significantly lower amount of parameters. Notably, the sampling time at $0.14$s is well within the sampling time of the system at $T_s = 5$ s as can be seen in Table \ref{tab:FourTank_Results} and is significantly faster than the sampling time of LSTM-MPC at $0.65$s.  

\begin{table}[h!]
\centering
\caption{Parameters and Performance for Mamba-MPC on the Quanser Aero2 Setup.}
\label{tab:AERO}
\begin{tabular}{l c}
\toprule
& \textbf{Quanser Aero2} \\
\cmidrule(lr){2-2}
& Mamba \\
\midrule
\multicolumn{2}{l}{\textit{Model Architecture}} \\
Model Dim. $D$            & 8 \\
Kernel Size               & 20 \\
State Dim.                & 4 \\
Horizon $N$               & 18 \\
Layers                    & 2 \\
Params                    & 983 \\
\midrule
\multicolumn{2}{l}{\textit{Controller Performance}} \\
MAE                       & $[87, 128, 134]\cdot10^{-3}$ \\ 
MSE                       & $[22, 90, 91]\cdot10^{-3}$ \\
Comp. Time (s)            & $0.038$ \\
Sampling Time (s)         & $0.04$ \\
\midrule
\multicolumn{2}{l}{\textit{Training Metrics}} \\
Time (min)                & 400 \\
Train Loss                & $2.2 \cdot 10^{-3}$ \\
Val. Loss                 & $2.0 \cdot 10^{-2}$ \\
\bottomrule
\end{tabular}
\end{table}

\begin{figure}
    \centering
    \includegraphics[width=1\linewidth]{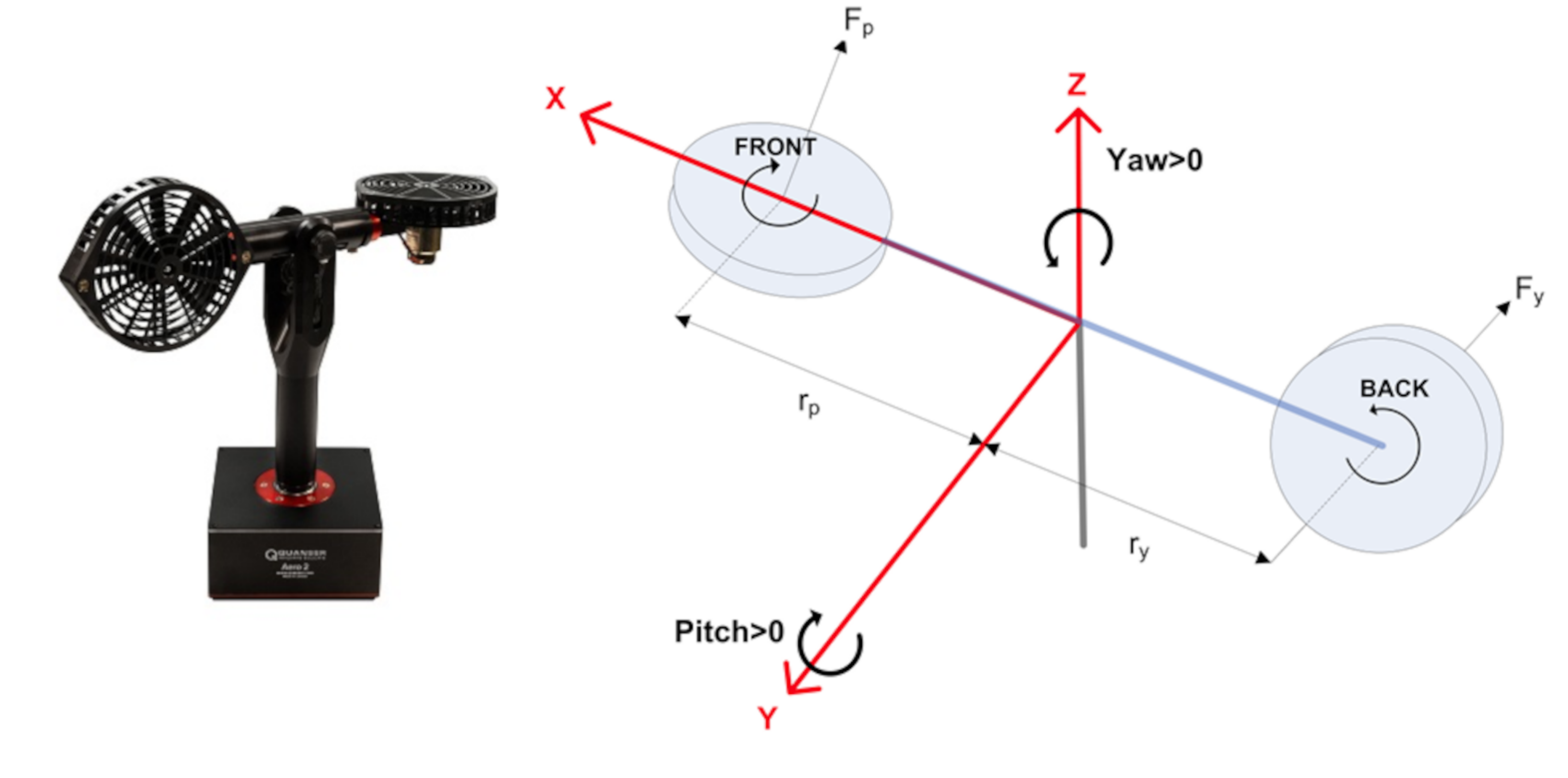}
    \caption{Left: Image of Quanser Aero2. Right: Free-body diagram of the dynamics of the Quanser Aero2}
    \label{fig:Aero2}
\end{figure}

\begin{figure}
    \centering
    \includegraphics[width=1\linewidth]{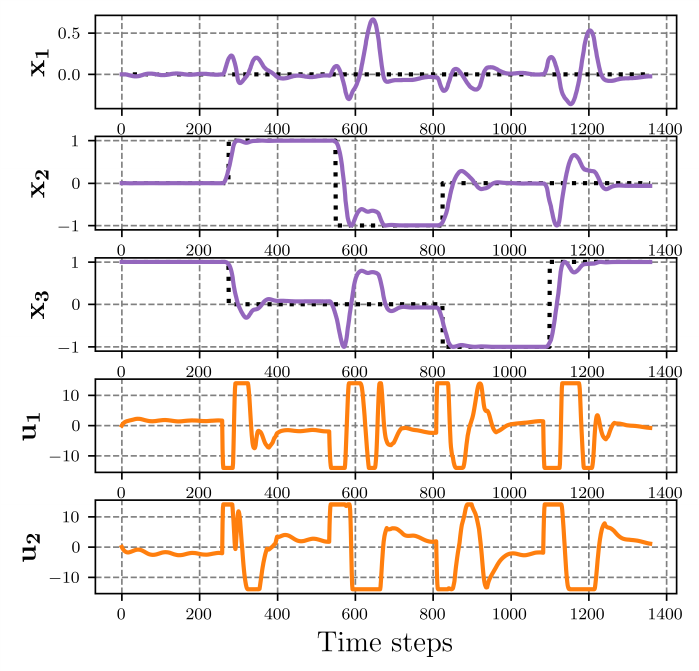}
    \caption{Closed-loop reference tracking of Mamba-MPC on the Quanser Aero2 physical setup. The states are indicated by (\textcolor{mypurple}{\rule[0.5ex]{1em}{1pt}}) which track a piece-wise constant reference indicated by (\textcolor{myblack}{\rule[0.5ex]{0.2em}{1pt}} \textcolor{myblack}{\rule[0.5ex]{0.2em}{1pt}}). The inputs are indicated by (\textcolor{myorange}{\rule[0.5ex]{1em}{1pt}}).}
    \label{fig:Aero2RefTrack}
\end{figure}
\section{Validation on a Physical Setup: The Quanser Aero2}\label{sec:AERO}
To assess the performance of Mamba-MPC on a physical setup, the Quanser Aero2 is considered, which can be seen in Figure \ref{fig:Aero2} alongside its free-body diagram. In this example, the 2-DOF configuration is considered. 
The Aero2 has inputs $u\in \mathbb{R}^{2}$ corresponding to the voltages supplied to the DC motors of rotor 1 and rotor 2. The system outputs $y\in \mathbb{R}^{3}$ encompassing the pitch angle, sine and cosine components of the yaw. For the initial conditions, pitch rate and yaw rate are added such that the initial condition $x_0\in \mathbb{R}^{5}$.
A dedicated training dataset is constructed via an open-loop experiment. Data acquisition is performed at a sampling frequency of $f_s = 100$ Hz for $T=240\cdot10^3$ samples where the input signal is constructed with a multi-level \emph{Pseudo Random Binary Signal} (PRBS) as carrier wave,  with amplitudes $[\pm 8, \pm 10, \pm 12]$ and maximum switching frequency of $2$ Hz. On top of the carrier wave, a multisine is superimposed featuring amplitudes $4,2,1$, respectively, across a frequency range of $0.015$ Hz to $2.485$ Hz containing $500$ harmonics. Due to computational limitations, the data was down-sampled to $f_s=25$ Hz to reduce the computational requirement. The definitive dataset is constructed according to \eqref{eq:StateEmbed} for a prediction horizon of $N=18$. Mamba-MPC is trained for $2000$ epochs according to \eqref{eq:LearningProblem} with hyperparameters given in Table \ref{tab:AERO}.

The model is implemented into the MPC scheme \eqref{prob:Mamba-MPC} with $Q=\text{diag}\left(6,1.5,1.5\right)$, $R=\text{diag}\left(0.1, 0.1\right)$ and $P=1\cdot 10^3$ with input constraints $\mathbb{U} = \{u\in \mathbb{R}^2| \text{col}(-14,-14) \leq u \leq \text{col}(14,14)\}$ and pitch angle constraints $\mathbb{Y}_1 = \{x_1 \in \mathbb{R}| -0.7 \leq x_1 \leq 0.7\}$. To remain within the sampling frequency of $f_s=25$Hz an SQP solver is used in CasADI. To further improve the computational speed, a warm start of the Lagrangian multipliers is added. Consider the closed-loop response for a piece-wise constant reference signal in Figure \ref{fig:Aero2RefTrack}. In general, Mamba-MPC is able to successfully track the reference signal for the pitch angle and yaw angle. In particular, the mean computation time per iteration is $0.038$ s which is below the sampling time of the system at $0.04$ s. Overall, Mamba-MPC is able to successfully stabilize and track a piece-wise constant reference on the Aero2 while remaining computationally tractable.\footnote{A video demonstration of Mamba-MPC on the Quanser Aero2 is available at: \url{https://youtu.be/JZUzMLH1HeU}}

\section{Conclusion}\label{sec:Conclusion}
In this paper, we developed Mamba-MPC, an MPC formulation that utilizes Mamba NN as a multi-step predictor, and a novel data embedding for Seq2Seq modeling. A complete and consistent mathematical description of the Mamba NN was given and the modeling advantages of Mamba were highlighted. We provided an embedding of initial conditions and input sequence such that decoder-only architectures can be utilized for multi-step prediction. We evaluated the performance of Mamba-MPC on two simulated nonlinear systems where we tested stabilization, noise robustness, reference tracking and MIMO reference tracking. Furthermore, we implemented Mamba-MPC on a physical Quanser Aero2 setup. On the simulated benchmarks, Mamba-MPC successfully stabilizes a nonlinear systems, it is robust to noisy data, outperforms LSTM-MPC on reference tracking, and scales favorably in comparison to LSTM-MPC and Transformer-MPC. Furthermore, on all benchmarks, Mamba-MPC stays within the desired computation time. Moreover, Mamba-MPC successfully performed MIMO reference tracking on the Quanser Aero2 at real-time control sampling interval of $0.04$s. Future work will consider:
\begin{itemize}
    \item Evaluating the performance of the SSM function when the diagonal terms in $A$ are allowed to be complex. Complex valued diagonal terms could increase the modeling capacity of Mamba in dynamical systems;
    \item Implementing Mamba-MPC using \texttt{Acados} for further computational efficiency;
    \item Developing a closed-loop stability analysis framework for Mamba-MPC.
\end{itemize}



\bibliographystyle{IEEEtran}
\bibliography{References}

\begin{IEEEbiography}[{\includegraphics[width=\linewidth,height=1.25in,clip,keepaspectratio,trim=0pt 0pt 0pt 0pt]{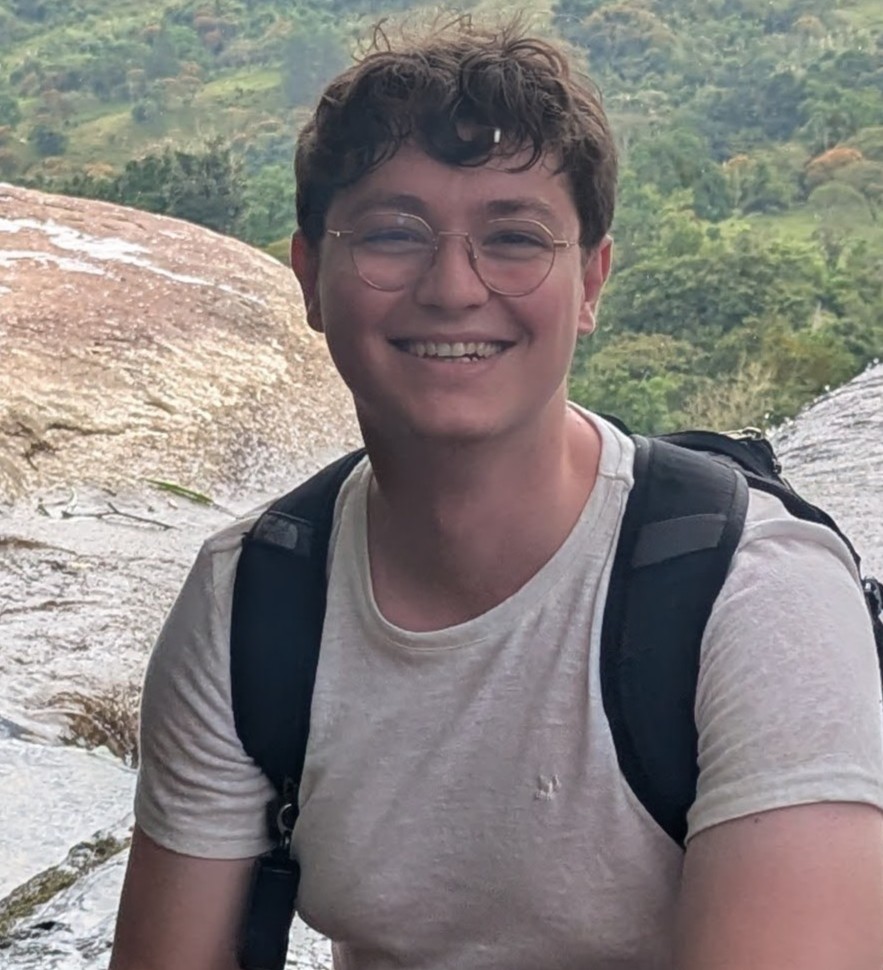}}]{Michiel Cevaal}{\space}
received the B.S. degree in Electrical Engineering from Eindhoven University of Technology,  Eindhoven, The Netherlands, in 2022 and the M.S. degree in Electrical Engineering from Eindhoven University of Technology, Eindhoven, The Netherlands in 2025. He is currently pursuing the Ph.D. degree at Tilburg University in collaboration with the Data Science Center of Excellence of the Ministry of Defense. His research is focused on autonomous multi-agent drone swarms in adversarial conditions. 
\end{IEEEbiography}

\begin{IEEEbiography}[{\includegraphics[width=\linewidth,height=1.25in,clip,keepaspectratio,trim=0pt 0pt 0pt 0pt]{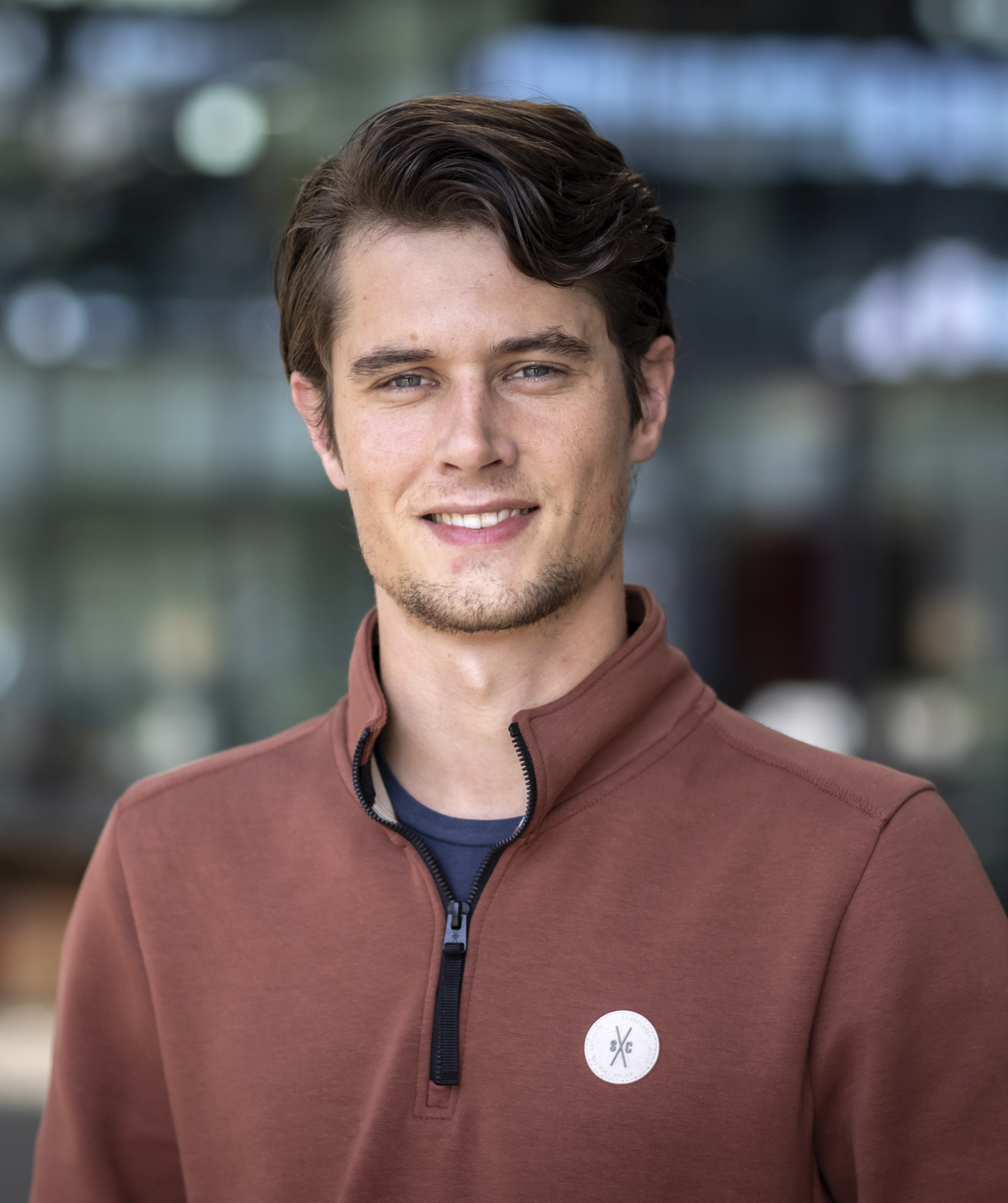}}]{Thomas O. de Jong}{\space}(Student Member, IEEE) 
received the B.S. degree in Mechanical Engineering from Eindhoven University of Technology, Eindhoven, The Netherlands, in 2021 and the M.S. degree in Systems and Control from Eindhoven University of Technology, Eindhoven, The Netherlands, in 2023. He is currently pursuing the Ph.D. degree in Systems and Control in Eindhoven, The Netherlands, in 2021 and the M.S. degree in Systems and Control at Eindhoven University of Technology, Eindhoven, The Netherlands.
\end{IEEEbiography}

\begin{IEEEbiography}[{\includegraphics[width=\linewidth,height=1.25in,clip,keepaspectratio,trim=0pt 0pt 0pt 0pt]{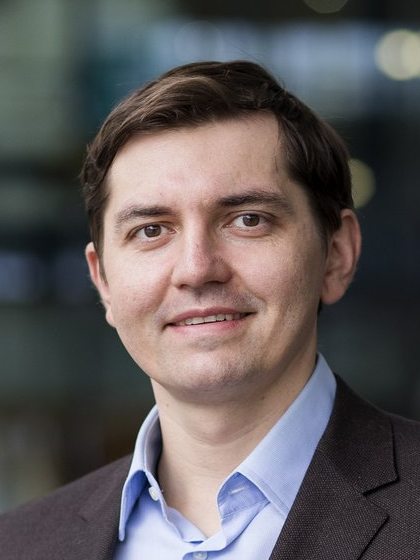}}]{Dr. Mircea Lazar}{\space}(Senior Member, IEEE)  is a Full Professor in Constrained Control of Complex Systems at the Electrical Engineering Department, Eindhoven University of Technology, The Netherlands. His research interests cover Lyapunov functions, distributed control, neural networks and constrained control of nonlinear and hybrid systems, including model predictive control. His research is driven by control problems in power systems, power electronics, high-precision mechatronics, automotive and biological systems. Dr. Lazar received the European Embedded Control Institute Ph.D. Award in 2007 for his Ph.D. dissertation. He received a VENI personal grant from the Dutch Research Council (NWO) in 2008 and he supervised 13 Ph.D. researchers (2 received the Cum Laude distinction). Dr. Lazar chaired the 4th IFAC Conference on Nonlinear Model Predictive Control in 2012. He is an Active Member of the IFAC Technical Committees 1.3 Discrete Event and Hybrid Systems, 2.3 Nonlinear Control Systems and an Associate Editor of IEEE Transactions on Automatic Control.
\end{IEEEbiography}

\end{document}